%% file: root.tex
\documentclass[letterpaper, 10 pt, journal]{IEEEtran}  % Comment this line out if you need a4paper

\IEEEoverridecommandlockouts                              % This command is only needed if 
\input{macros.tex}

\title{\LARGE \bf
	Routing for Traffic Networks with Mixed Autonomy
}

\author{Daniel A. Lazar$^{1}$, Samuel Coogan$^{2}$, and Ramtin Pedarsani$^{3}$% <-this % stops a space
	\thanks{$^{1}$Daniel Lazar is with the Department of Electrical and Computer Engineering, 
		University of California, Santa Barbara
		{\tt\small dlazar@ece.ucsb.edu}}%
	\thanks{$^{2}$Samuel Coogan is with the School of Electrical and Computer Engineering and the School of Civil and Environmental Engineering, 
		Georgia Institute of Technology
		{\tt\small sam.coogan@gatech.edu}}%
	\thanks{$^{3}$Ramtin Pedarsani is with the Department of Electrical and Computer Engineering, 
		University of California, Santa Barbara
		{\tt\small ramtin@ece.ucsb.edu}}%
}

\begin{document}
	
\maketitle
\thispagestyle{empty}
\pagestyle{empty}

%%%%%%%%%%%%%%%%%%%%%%%%%%%%%%%%%%%%%%%%%%%%%%%%%%%%%%%%%%%%%%%%%%%%%%%%%%%%%%%%
\begin{abstract}
	
	In this work we propose a macroscopic model for studying routing on networks shared between human-driven and autonomous vehicles that captures the effects of autonomous vehicles forming platoons. We use this to study inefficiency due to selfish routing and bound the Price of Anarchy (PoA), the maximum ratio between total delay experienced by selfish users and the minimum possible total delay. To do so, we establish two road capacity models, each corresponding to an assumption regarding the platooning capabilities of autonomous vehicles. Using these we develop a class of road delay functions, parameterized by the road capacity, that are polynomial with respect to vehicle flow. We then bound the PoA and the bicriteria, another measure of the inefficiency due to selfish routing. We find these bounds depend on: 1) the degree of the polynomial in the road cost function and 2) the degree of asymmetry, the difference in how human-driven and autonomous traffic affect congestion. We demonstrate that these bounds recover the classical bounds when no asymmetry exists. We show the bounds are tight in certain cases and that the PoA bound is order-optimal with respect to the degree of asymmetry.
	
\end{abstract}

%%%%%%%%%%%%%%%%%%%%%%%%%%%%%%%%%%%%%%%%%%%%%%%%%%%%%%%%%%%%%%%%%%%%%%%%%%%%%%%%
\section{INTRODUCTION}
\input{./intro.tex}

\section{RELATED WORK}\label{sct:prev_works}
\input{./related_work.tex}

\section{NETWORK MODEL}\label{sct:modeling}
\input{./network_model.tex}

\section{BOUNDING THE PRICE OF ANARCHY}\label{sct:main_results}
\input{./bounds.tex}

\section{ESTABLISHING LOWER BOUNDS BY EXAMPLE}\label{sct:tightness}
\input{./discussion_examples.tex}

\section{CONCLUSION}\label{sec:conclusion}
\input{./conclusion.tex}

\input{root.bbl}
\section{APPENDIX}
\input{./proofs.tex}

\input{./bios.tex}

\addtolength{\textheight}{-3cm}   % This command serves to balance the column lengths

%%%%%%%%%%%%%%%%%%%%%%%%%%%%%%%%%%%%%%%%%%%%%%%%%%%%%%%%%%%%%%%%%%%%%%%%%%%%%%%%

\end{document}

%% file: macros.tex
\usepackage{subcaption}
\usepackage{graphicx}
\usepackage{enumerate}
\usepackage{float}
\usepackage[font=footnotesize]{caption}

\usepackage{pgfplots}
\usepackage{amsmath,amssymb}
\usepackage{amsthm}

\usepackage{tikz}
\usepackage{comment}
\usepackage{cite}

\usepackage{mathtools, cuted}
\usepackage{multicol}

\usetikzlibrary{automata,arrows,positioning,calc}

\usepackage{todonotes}

\newtheorem{theorem}{Theorem}
\newtheorem{corollary}{Corollary}
\newtheorem{lemma}{Lemma}
\newtheorem{remark}{Remark}
\newtheorem{example}{Example}
\newtheorem{definition}{Definition}
\newtheorem{proposition}{Proposition}

\newcommand{\der}{\mathrm{d}}
\newcommand{\xeq}{x^{\text{EQ}}}
\newcommand{\yeq}{y^{\text{EQ}}}
\newcommand{\zeq}{z^{\text{EQ}}}
\newcommand{\qeq}{q^{\text{EQ}}}
\newcommand{\xopt}{x^{\text{OPT}}}
\newcommand{\yopt}{y^{\text{OPT}}}
\newcommand{\zopt}{z^{\text{OPT}}}
\newcommand{\qopt}{q^{\text{OPT}}}

\newcommand{\cagg}{c^{\text{AGG}}}
\newcommand{\Cagg}{C^{\text{AGG}}}

\newcommand{\feq}{f^{\text{EQ}}}
\newcommand{\fopt}{f^{\text{OPT}}}

\newcommand{\nidcs}{{[}n{]}}
\newcommand{\ffli}{t_i} % free-flow latency on road i
\newcommand{\shdwy}{\bar{h}} % short headway
\newcommand{\shdwyi}{\bar{h}_i} % short headway on road i
\newcommand{\lhdwy}{h} % long headway
\newcommand{\lhdwyi}{h_i} % long headway on road i

\newcommand{\costscli}{\rho_i}
\newcommand{\feasrout}{\mathcal{X}}

\DeclareMathOperator{\E}{\mathbb{E}}

%% file: intro.tex
In recent years, automobiles are increasingly equipped with autonomous and semi-autonomous technology, which has potential to dramatically decrease traffic congestion \cite{dot:2015zr}. Specifically, autonomous technologies enable \emph{platooning}, in which these vehicles automatically maintain short headways between them via adaptive cruise control (ACC) or cooperative adaptive cruise control (CACC). ACC uses sensing such as radar or LIDAR to maintain a specific distance to the preceeding vehicle with faster-than-human reaction time, and CACC augments this with vehicle-to-vehicle communications.
	
When all vehicles are autonomous, the use of platooning has the potential to increase network capacity as much as three-fold \cite{lioris2017platoons}, smooth traffic flow by eliminating shockwaves of slowing vehicles \cite{stern2018dissipation}, and enable synchronous acceleration at green lights\cite{Askari:2016fy}. However, the presence of human-driven vehicles -- leading to \emph{mixed autonomy} -- makes much of these benefits unclear.

\begin{figure}
	\centering
	\begin{subfigure}[b]{1\linewidth}
		\includegraphics[width=1\linewidth]{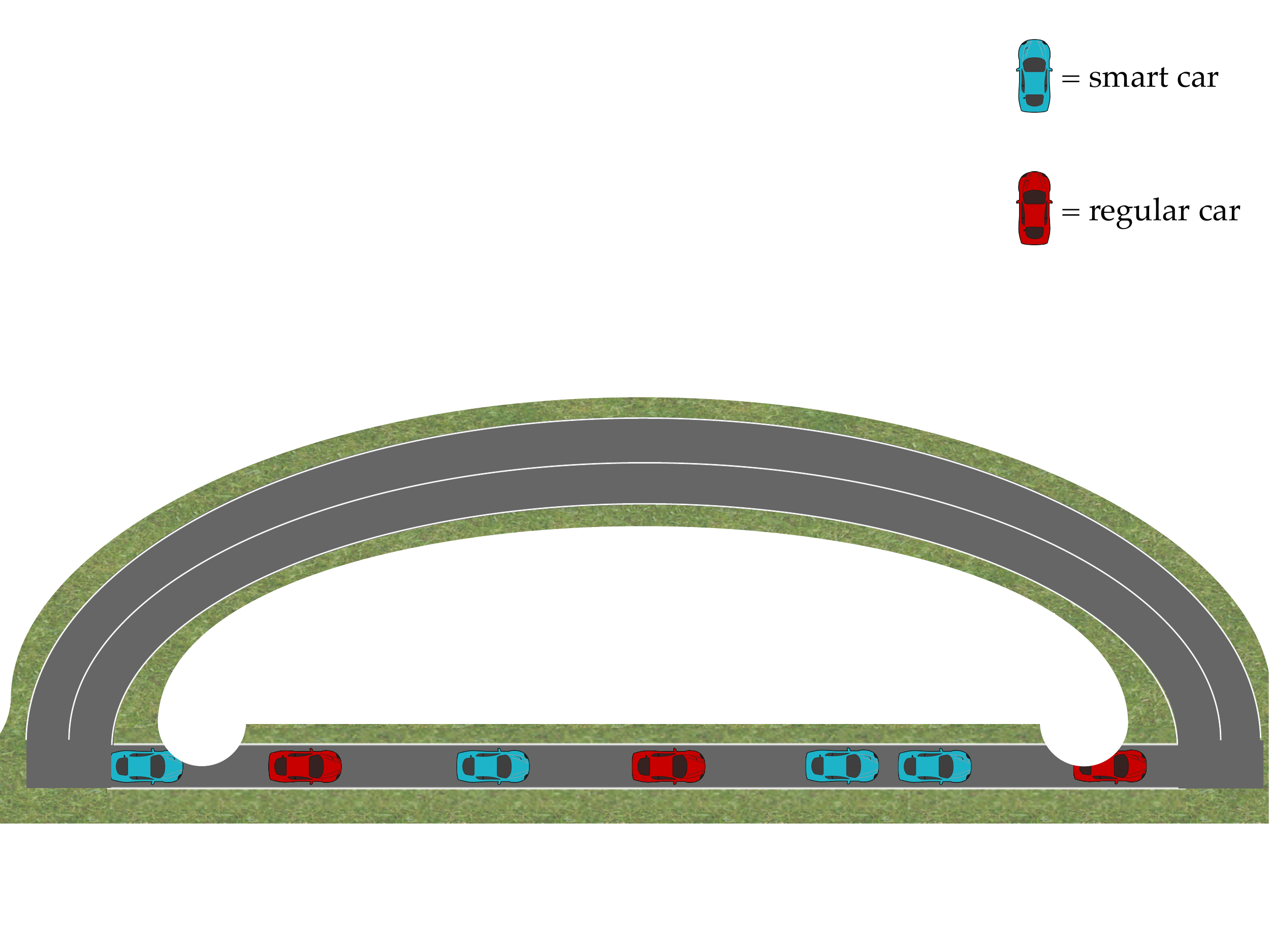}
		\caption{}
	\end{subfigure}
	\begin{subfigure}[b]{1\linewidth}
		\includegraphics[width=1\linewidth]{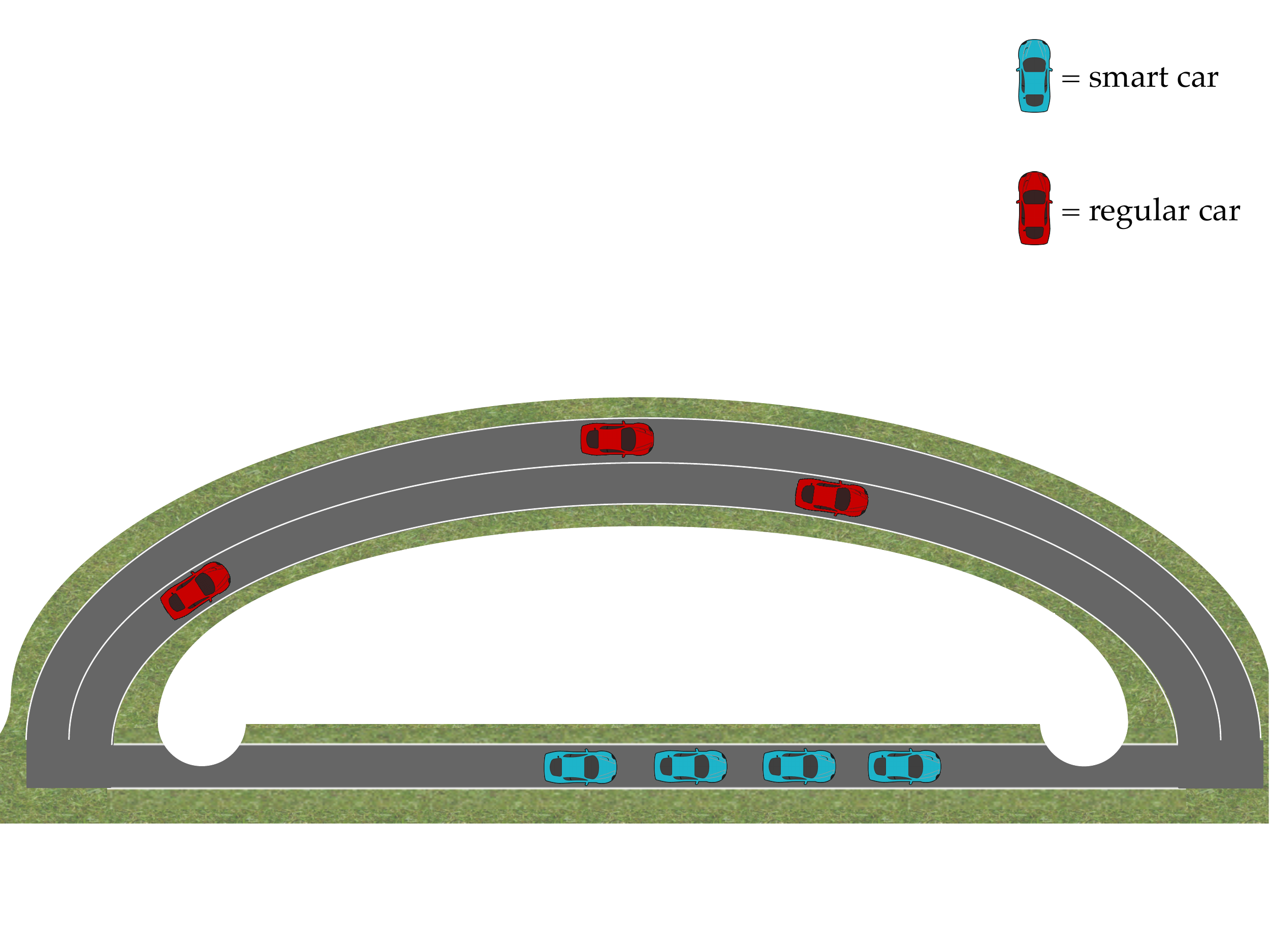}
		\caption{}
	\end{subfigure}
	\caption{A social planner can decrease overall travel times by make routing decisions that utilize autonomous vehicles' ability to platoon, and choosing different routes for human-driven vehicles (red) and autonomous vehicles (blue). When vehicles route selfishly (a), vehicles pack onto a congested road. In optimal routing (b), only autonomous vehicles are sent onto the road most amenable to platooning.}
	\label{fig:frontfig}
\end{figure}

Moreover, even in the absence of autonomous capabilities, it is well known that if drivers route selfishly and minimize their individual traffic delays, this does not in general minimize \emph{overall} traffic delay. Understanding the extent of this phenomenon can help city planners -- if selfish routing does not adversely affect congestion too much, then it may not be necessary to try to control vehicle flow using schemes such as tolling. Alternatively, if selfishness can lead to much worse congestion, then a city planner may wish to try to control human routing decisions. See Figure~\ref{fig:frontfig} for an example of selfish routing and optimal routing in mixed autonomy.

The ratio between traffic delay under worst-case selfish routing and optimal routing is called the Price of Anarchy (PoA) and is well understood for networks with only human-driven vehicles \cite{roughgarden2002bad, roughgarden2002selfishthesis, roughgarden2003topology, roughgarden2005selfish, correa2004selfish, correa2008geometric}. Many such works also bound the bicriteria, which quantifies, for any given volume of vehicle flow demand, how much additional flow can be routed optimally to result in the same overall latency as the original volume of traffic routed selfishly. Other studies have bounded the PoA with multiple modes of transportation \cite{perakis2007price, chau2003price}. However, these prior works require assumptions that do not capture the effects of congestion on roads shared between human-driven and autonomous vehicles, leaving open the question of the price of anarchy in mixed autonomy. In fact, \emph{we show that these previous results do not hold, and the PoA for roads with mixed autonomy can in general be unbounded!}

Motivated by this observation, in this paper we provide novel bounds on the PoA and bicriteria that depend on the degree to which platooning decreases congestion, as well as the degree of the polynomial describing road delay. To do so, we use two models that describe road capacity as a function of the fraction of vehicles on the road that are autonomous; each model corresponds to a different assumption regarding the technology that enables platooning. We use these capacity models with a known polynomial road delay function, and, for this class of latency functions, we bound the price of anarchy and bicriteria. We develop two mechanisms for bounding the PoA, which yield bounds that are tighter depending on platoon spacing and polynomial degree. In our development we provide the main elements of our proofs and defer proofs of lemmas to the appendix. 

In our formulation, the benefit due to the presence of autonomous vehicles is limited to platoon formation, and the probability that each vehicle is autonomous is independent of the surrounding vehicles. While we acknowledge that autonomous vehicles yield other benefits such as smoothing traffic shockwaves, we consider platooning because it is a mature technology that is commercially available. Further, if autonomous vehicles actively rearrange themselves to form platoons, the resulting capacity falls between the two capacity models presented here \cite{lazar2018influence}. \\

\begin{figure}
	\begin{center}
		\begin{tikzpicture}[->, >=stealth', auto, semithick, node distance=4cm]
		\tikzstyle{every state}=[fill=white,draw=black,thick,text=black,scale=1]
		\node[state]    (0)               {$s$};
		\node[state]    (1)[right of=0]   {$t$};
		\path
		(0) edge[bend left]		node{$c_1(x,y) = 1$}     (1)
		(0) edge[bend right]	node[below]{$c_2(x,y) = \zeta x$}     (1);
		\end{tikzpicture}
	\end{center} 
	\caption{Road network with price of anarchy and bicriteria that grow unboundedly with $\zeta$ when considering $1/\zeta$ units of human-driven flow and $1$ unit of autonomous flow demand, with $\zeta \ge 1$. Function arguments $x$ and $y$ respectively denote human-driven and autonomous flow on a road.}
	\label{fig:unbouded2road}
\end{figure}
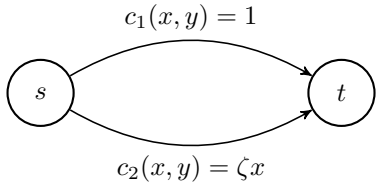

\noindent\textbf{Motivating Example.} To show that bounds for the price of anarchy previously developed for roads with only one type of vehicle (\emph{i.e.} no autonomous vehicles) do not hold, we present an example of a road network with unbounded price of anarchy (Fig.~\ref{fig:unbouded2road}). Consider a network of two parallel roads, with road latency functions $c_1(x,y)=1$ and $c_2(x,y)=\zeta x$. On each road the latency is a function of the human-driven flow ($x$) and the autonomous flow ($y$) on that road. Suppose we have $\frac{1}{\zeta}$ units of human-driven vehicle flow and $1$ unit of autonomous traffic demand to cross from node $s$ to node $t$, with $\zeta \ge 1$.

Optimal routing puts all human-driven cars on the top road and all autonomous cars on the bottom road; when vehicles route selfishly they all end up on the bottom road. This yields a price of anarchy of $\zeta+1$. The bicriteria is also $\zeta + 1$, as $\zeta+1$ times as much traffic, optimally routed, yields the same total cost as the original amount of traffic at Wardrop Equilibrium.\footnote{Though in this case autonomous vehicles do not affect congestion, other examples in Section \ref{sct:tightness} also yield an unbounded PoA with both vehicle types affecting congestion.} This examples leads us already to our first proposition, which lays the foundations for the contributions of this paper:
\begin{proposition}
	The price of anarchy and bicriteria are in general unbounded in mixed autonomy.
\end{proposition}

Motivated by this proposition, we develop the notion of the \emph{degree of asymmetry} of a road and use this, in conjunction with the degree of the polynomial cost function, to parameterize the bound on the price of anarchy. To summarize, we
\begin{enumerate}
	\item show that previous PoA results do not hold for mixed autonomy,
	\item develop a realistic class of polynomial cost functions for traffic of mixed autonomy,
	\item develop two mechanisms for bounding the PoA using this cost function with both capacity models, and
	\item bound the price of anarchy and bicriteria and analyze the tightness of our bounds.
\end{enumerate}

Some of the contributions above are related to our previous work. In one previous paper, we use similar capacity models along with a latency function based on M/M/1 queues to find optimal routing for human-driven and autonomous vehicles over two parallel roads \cite{lazar2017routing}. In another work, we consider maximizing capacity, using the second capacity model in this paper, via a sequence of vehicle reorderings in which autonomous vehicles influence human drivers \cite{lazar2018influence}. 

In this paper we expand on conference version, which considers the price of anarchy and bicriteria in mixed autonomy with affine latency functions \cite{lazar2018poa}. The derivation considers the first capacity model in conjunction with a latency function that can be considered as a special case of the latency function presented in this paper, with polynomial degree one. The bound presented in that paper corresponds to that of Theorems \ref{thm:low_asym} and \ref{thm:bicriteria}, limited to polynomials of degree one.  In the current work we consider a broader class of latency functions that includes both capacity models and latency functions with arbitrary polynomial degree. We expand the previous bounds to include this broader class of latency functions. However, this bound only holds under certain regimes for the allowable parameters of the capacity model and latency functions. In response, we present a PoA bound that holds in all regimes and is tighter even in some regimes in which the previous bound holds.

%% file: related_work.tex
\subsection{Congestion Games and Wardrop Equilibria}

Our work is related to the optimal traffic assignment problem, \emph{e.g.} \cite{dafermos1969traffic_general}, which studies how to optimally route vehicles on a network when the cost (\emph{i.e.} delay) on a road link is a function of the vehicles that travel on that link. We are concerned specifically with the relationship between optimal traffic assignment and Wardrop Equilibria, which occur when drivers choose their paths selfishly. For a survey on literature on Wardrop Equilibria, see \cite{correa2011wardrop}; \cite{depalma1998optimization} describes other notions of equilibria. Classic works on Wardrop Equilibria and the associated tools for analyzing them include \cite{smith1979existence, dafermos1980variationalinequalities, florian1995network}.

In an important development, Smith \cite{smith1979existence} establishes the widely used \emph{variational inequality} and uses it to describe flows at Wardrop Equilibrium, in which all users sharing an origin and destination use paths of equal cost and no unused path has a smaller cost. For any feasible flow $z$ and equilibrium flow $\zeq$, the variational inequality dictates that
\begin{align}\label{eq:VI}
\langle c(\zeq), \zeq-z \rangle \le 0 \; ,
\end{align}
where $z$ is a vector describing vehicle flow on each road and $c(z)$ maps a vector of flows to a vector of the delay on each road. The Variational Inequality is fundamental for establishing our PoA bound. In the case of mixed autonomy each road is duplicated to account for the two types of traffic; a full description is provided in Section \ref{sct:modeling}.

\subsection{Multiclass Traffic}

Some previous works consider traffic assignment and Wardrop Equilibria with multiclass traffic, meaning traffic with multiple vehicle types and transportation modes that affect and experience congestion differently (\emph{e.g.} \cite{dafermos1972multiclass_user, hearn1984convex, florian1977traffic}). 

Florian \cite{florian1977traffic} demonstrates how to calculate equilibria for a multimodal system involving personal automobiles and public transportation. They use a relaxation that assumes public transportation will take the path that would be shortest in the absence of cars. In the case of mixed autonomy, this is not a fair assumption. 

Dafermos \cite{dafermos1972multiclass_user} assumes the Jacobian of the cost function is symmetric and positive definite. Similarly, Hearn \emph{et. al.} \cite{hearn1984convex} deal with a monotone cost function, \emph{i.e.} satisfying the property
\begin{align}\label{eq:monotonicity}
	\langle c(z)-c(q),z-q \rangle \ge 0 
\end{align}
for flow vectors $z$ and $q$.

However, traffic networks with mixed autonomy are in general nonmonotone. To see this, consider a network of two roads with costs $c_1(x,y)=3x +y + t_1$ and $c_2(x,y)=3x+2y+t_2$, where $t_1$ and $t_2$ are constants denoting the free-flow latency on roads 1 and 2. This corresponds to a road in which autonomous vehicles can platoon closely and another road on which they cannot platoon as closely. The Jacobian of the cost function is as follows:
\begin{align*}
	\begin{bmatrix}
	3 & 1 & 0 & 0 \\
	3 & 1 & 0 & 0 \\
	0 & 0 & 3 & 2 \\
	0 & 0 & 3 & 2
	\end{bmatrix} \; ,
\end{align*}
which is not symmetric, and the vector $z = \begin{bmatrix}
-1 & 2 & 0 & 0 \end{bmatrix}^T$ demonstrates that it is also not positive semidefinite. Monotonicity is closely related to the positive (semi-) definiteness of the Jacobian of the cost function. To show that the monotonicity condition is violated as well, consider that there are 2 units of human-driven flow demand and 3 units autonomous flow demand. With one routing in which all human-driven flow is on the first road and all autonomous flow is on the second and another routing with these reversed, we find that the monotonicity condition is violated. 

Similarly, Faroukhi \emph{et. al.} \cite{farokhi2014potential} prove that in heterogeneous routing games with cost functions that are continuously differentiable, nonnegative for feasible flows, and nondecreasing in each of their arguments, then at least one equilibrium is guaranteed to exist. These mild conditions are satisfied in our setting. For heterogeneous games with two types, they further prove a necessary and sufficient condition for a potential function (and therefore unique equilibrium) to exist. However, the condition required can be considered a relaxation of the condition that the Jacobian of the cost function be symmetric. While broader than strict symmetry, this condition is still not satisfied in mixed autonomy. Notably, they describe tolls that, when applied, yield a cost function that satisfies this condition.

As described above, these previous works in multiclass traffic require restrictive assumptions and therefore do not apply to the case of mixed autonomy. In fact, in the case of mixed autonomy, the routing game is not formally a proper congestion game, as it cannot be described with a potential function. Nonetheless, in this paper we adapt tools developed for such games to derive results for mixed autonomous traffic.

\subsection{Price of Anarchy}

There is an abundance of research into the price of anarchy in nonatomic congestion games, codified in \cite{roughgarden2002bad, roughgarden2002selfishthesis, roughgarden2003topology, roughgarden2005selfish, correa2004selfish, correa2008geometric}. In \cite{correa2008geometric}, the authors develop a general tool for analyzing price of anarchy in nonatomic congestion games. Though their development is specific to monotone cost functions, in this paper we broaden it to cost functions that are not necessarily monotone. Also relatedly, we find that in the case of no asymmetry, our price of anarchy bound for polynomial cost functions simplifies to the classic bound in \cite{roughgarden2002selfishthesis, roughgarden2003topology, roughgarden2005selfish}.

The previously mentioned works consider primarily single-type traffic. Perakis \cite{perakis2007price} considers PoA in multiclass traffic using nonseparable, asymmetric, nonlinear cost functions with inelastic demand. However, they restrict their analysis to the case that the Jacobian matrix of the cost function is positive semidefinite. Similarly, Chau and Sim \cite{chau2003price} consider the PoA for multiclass traffic with \emph{elastic} demand with symmetric cost functions and positive semidefinite Jacobian of the cost function. As demonstrated earlier, these assumptions are violated in the case of mixed autonomy.

\subsection{Autonomy}

In one line of research, autonomous vehicles are controlled to \emph{locally} improve traffic by smoothing out stop-and-go shockwaves in congested traffic
\cite{Shladover:1978lo, Darbha:1999dw,Yi:2006hb, Pueboobpaphan:2010qe, Orosz:2016hb, cui2017stabilizing, stern2018dissipation, wu2017emergent, wu2018stabilizing, li2018modeling}, optimally sending platooned vehicles through highway bottlenecks \cite{motie2016throughput, li2018modeling}, and simultaneously accelerating platooned vehicles at signalized intersections \cite{Askari:2016fy, lioris2017platoons}. Other papers investigate fuel savings attained using autonomous vehicles \cite{liang2016heavy, adler2016optimal,  turri2017cooperative, van2018fuel} or jointly controlling vehicles on a highway to localize and eliminate traffic disturbances \cite{sivaranjani2015localization}. Some works consider optimally routing and rebalancing a fleet of autonomous vehicles \cite{zhang2016routing}, though these generally consider a simpler model for congestion, in which capacitated roads have constant latency for flows below their capacity.

Some previous works have related models for road capacity and throughput under mixed autonomy, in particular \cite{Askari:2016fy, askari2017effect}, \cite{lazar2017routing, lazar2018poa, lazar2018influence}. In our previous work \cite{lazar2018poa}, we found the price of anarchy for affine cost functions with one capacity model. This is a special case of our work in this paper, which goes beyond the previous work by considering polynomial cost functions and incorporating both capacity models. Further, we introduce a novel mechanism for finding the PoA in mixed autonomy, leading to a tighter bound than the one previously found.

%% file: network_model.tex
Consider a congestion game on a network of $N$ roads, with nonatomic drivers (meaning each control an infinitesimally small unit of vehicle flow) traveling across $m$ origin-destination pairs, each pair associated with $\alpha_i$ units of human-driven vehicle flow demand and $\beta_i$ units of autonomous vehicle flow demand. We use ${[}n{]}=\{1,2,\ldots,n\}$ to denote the set of roads. We fully describe driver behavior on a network by using a vector of vehicle flows, which describes the volume of vehicles of each type that travels on each road. This vector has size equal to twice the number of roads and uses alternating entries to denote human-driven and autonomous vehicle flow on a road. We use $x_i$ and $y_i$ to refer to human-driven and autonomous flow on road $i$, respectively. Then the flow vector $z$ is as follows:
\begin{align*}
z=\begin{bmatrix}
x_1 & y_1 & x_2& y_2 & \ldots & x_n & y_n
\end{bmatrix}^T \; \in \mathbb{R}^{2N}_{\ge 0} \; .
\end{align*}

We refer to this flow vector as a \emph{routing} or a \emph{strategy}. We use $\feasrout \subseteq \mathbb{R}^{2N}_{\ge 0}$ to denote the set of feasible routings, meaning routings that route all flow demand from their origin nodes to their destination nodes while respecting conservation of flow in the network.

When needing to distinguish between two vectors, we use $v_i$ and $w_i$ in place of $x_i$ and $y_i$ and $q$ in place of $z$. We assume that human-driven and autonomous vehicles experience congestion identically. To capture the differing effects of each type of flow on a road's latency, we construct cost function $c(z): \mathbb{R}^{2N}_{\ge 0} \rightarrow \mathbb{R}^{2N}_{\ge 0}$ as follows:
\begin{align*}
	c(z)= \begin{bmatrix}
	c_1(x_1, y_1) \\
	c_1(x_1, y_1) \\
	c_2(x_2, y_2) \\
	c_2(x_2, y_2) \\
	\ldots \\
	c_n(x_n, y_n) \\
	c_n(x_n, y_n) 
	\end{bmatrix} \; ,
\end{align*}
where $c_i(x_i,y_i)$ is the latency on road $i$ when $x_i$ units of human-driven vehicles and $y_i$ units of autonomous vehicles use the road. The social cost, which is the \emph{aggregate delay} experienced by all users of the network, is then $C(z):=\langle c(z),z\rangle$. In the following sections we develop models for road capacity and road delay in order to construct cost functions.

\subsection{Capacity Models}\label{sct:cap_models}

\begin{figure}
	\centering
	\includegraphics[width=0.9\linewidth]{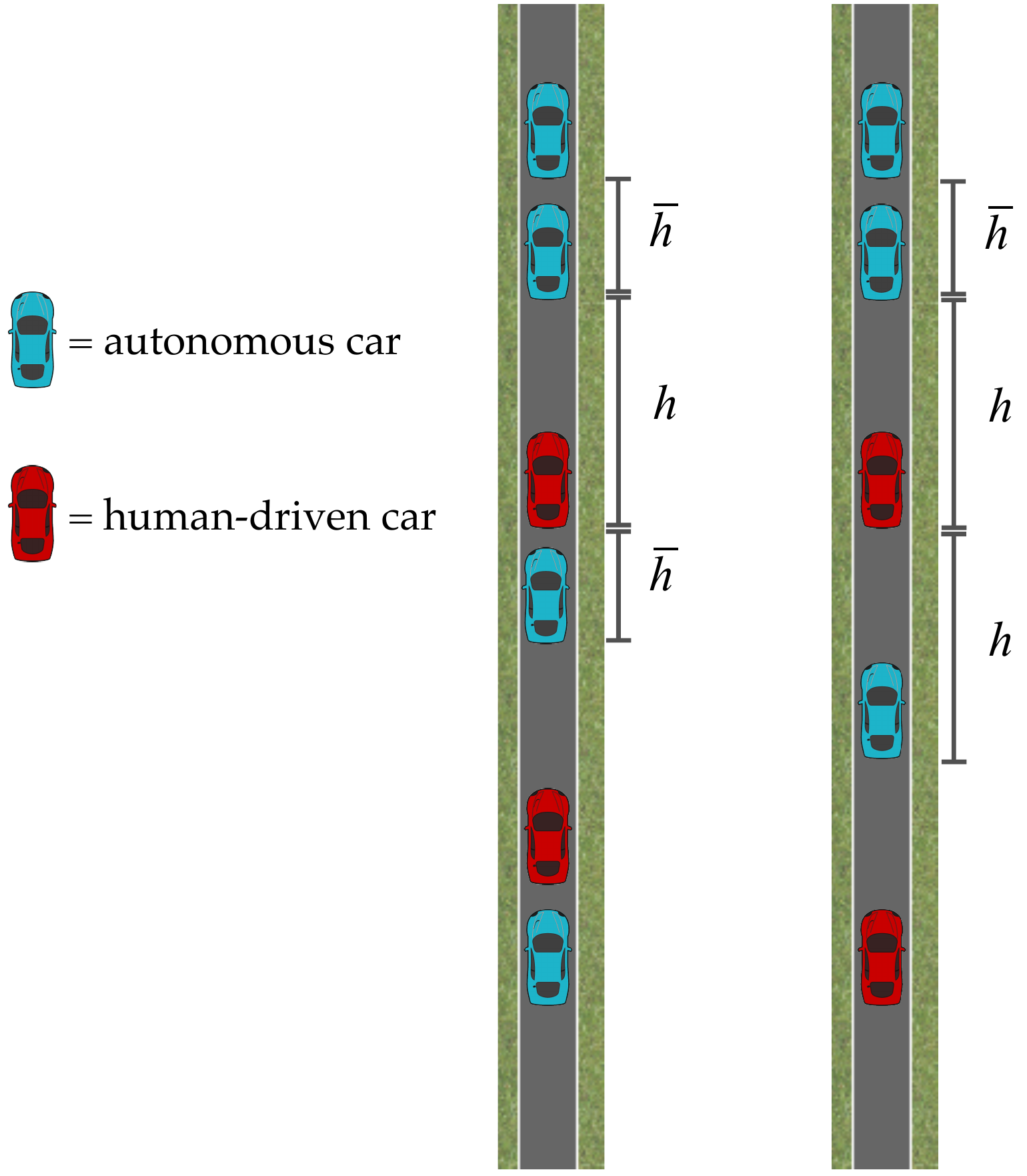}
	\caption{Capacity models 1 and 2. In capacity model 1 (left), autonomous cars can platoon behind any vehicle, and therefore take up length $\shdwy$. In capacity model 2 (right), autonomous vehicles can only platoon behind other autonomous vehicles; in that case they take up length $\shdwy$, but if following a human-driven vehicle, they take up length $\lhdwy$. Human-driven vehicles always take up length $\lhdwy$.}
	\label{fig:headways}
\end{figure}

We model the capacity of a road under two assumptions: 1) autonomous vehicles can platoon (follow closely) behind any vehicle and 2) autonomous vehicles can only platoon behind other autonomous vehicles. Let $d_i$ denote road length and let $\shdwyi$ and $h_i$ denote the nominal space taken up by a platooned and nonplatooned vehicle, respectively. The capacity will be a function of the \emph{autonomy level} of the road, denoted $\alpha(x_i,y_i)= y_i/(x_i+y_i)$. We define the capacity of a road as the number of vehicles, traveling at the road's nominal speed, that can fit on the length of the road. This is calculated by dividing the length of the road by the average space taken up by a car on the road, which is a function of autonomy level. We formalize this in the following proposition\footnote{Found in \cite{lazar2017routing}, contemporaneously in \cite{askari2017effect}.}.
\begin{proposition}
	Assume that vehicles are placed on a road as the result of a Bernoulli process. If autonomous vehicles can platoon behind any vehicle, therefore occupying road length $\shdwyi$ when traveling at nominal velocity, and human driven vehicles do not platoon (therefore occupying road length $\lhdwyi$ at nominal velocity), then the capacity is
	\begin{align}\label{eq:capm1}
		m_{i}(x_i,y_i) &= \frac{d_i}{\alpha(x_i,y_i) \shdwyi + (1-\alpha(x_i,y_i))\lhdwyi} \; .
	\end{align}
	If autonomous vehicles only platoon behind other autonomous vehicles and human driven vehicles cannot platoon, then the capacity is
	\begin{align}\label{eq:capm2}
		m_{i}(x_i,y_i) &= \frac{d_i}{\alpha^2(x_i,y_i)\shdwyi + (1-\alpha^2(x_i,y_i))\lhdwyi} \; .
	\end{align}
\end{proposition}
\begin{proof}
	We first justify the proposition for the first capacity model. Autonomous vehicles follow any vehicle with the same headway (occupying total space $\shdwyi$), as do human-driven vehicles with a different headway (occupying $\lhdwyi$). The space taken up by an average vehicle, as the number of vehicles grows large, is a weighted combination of those two spacings that depends on the autonomy level. Note that this capacity model does not depend on the ordering of the vehicles.
	
	For the second capacity model, we assume the vehicles are placed as the result of a Bernoulli process with parameter $\alpha_i$. Consider $M$ vehicles, each with length $L$, with $s_m$ denoting the headway of vehicle $m$. Note that the front vehicle will have $s_m=0$. The expected total space taken up is, due to linearity of expectation,
	\begin{align*}
	\E[\sum_{m=1}^{M} L &+ s_m] = M L + \sum_{m=1}^{M-1}\E[s_m] \\
	&=  (M-1)(\alpha^2(x_i,y_i)\shdwyi + (1-\alpha^2(x_i,y_i))\lhdwy) + L \; .
	\end{align*}
	Then, as the number of vehicles grows large, the average space occupied by a vehicle approaches $\alpha^2(x_i,y_i)\shdwyi + (1-\alpha^2(x_i,y_i))\lhdwy$, yielding the above expression for capacity model 2.
\end{proof}

Figure~\ref{fig:headways} provides an illustration of the technology assumptions.To make the meaning of nominal vehicle spacing more concrete, we offer one way of calculating spacing: let $L$ denote vehicle length, $\tau_{h,i}$ and $\tau_{a,i}$ denote the reaction speed of human-driven and autonomous vehicles, respectively. Let $v_i$ be the nominal speed on road $i$, which is likely the road's speed limit. Then we consider $\shdwy_i = L+v_i \tau_{a,i}$ and $\lhdwy_i = L+ v_i \tau_{h,i}$.\footnote{In general we consider $h_i\ge \bar{h}_i$, but we do not formally make this assumption. Our theoretical results hold even in the case that $\lhdwy_i \ge \shdwy_i$ on some roads and $\lhdwy_j < \shdwy_j$ on others.}

\subsection{Delay Model}

We now propose a model for the delay incurred by mixed traffic resulting from the capacity models derived in the previous subsection. We base our road delay function on that used in the Traffic Assignment Manual, published in 1964 by the Bureau of Public Roads (BPR) \cite{manual1964bureau, florian1977traffic, sheffi1985urban}. In a review of link capacity and cost functions \cite{branston1976link}, this cost function is analyzed and compared to alternate analytically difficult cost functions; the authors find the BPR cost function incurs little or no loss in accuracy compared to more complicated models.

We are concerned with the total, or aggregate, delay for all drivers on the network. Motivated by the delay function described in the Bureau of Public Roads' Traffic Assignment Manual \cite{manual1964bureau, florian1977traffic, sheffi1985urban}, we propose the following road delay function for mixed autonomous traffic:
\begin{align}\label{eq:cost_fn}
c_i(x_i,y_i) = \ffli(1+\costscli(\frac{x_i+y_i}{m_i(x_i,y_i)})^{\sigma_i}) \; .
\end{align}

Here, $\ffli$ denotes the free-flow delay on road $i$, and $\costscli$ and $\sigma_i$ are model parameters. Typical values for $\costscli$ and $\sigma_i$ are 0.15 and 4, respectively \cite{sheffi1985urban}. However, our solution methodology is valid for any parameters such that $\ffli \ge 0$, $\costscli \ge 0$ and $\sigma_i \ge 1$. We extend the BPR cost function by considering total flow on a road, and parameterize the capacity by the autonomy level of the flow on a road.

\begin{remark}
	This choice of cost functions implies that road delay is separable, meaning that the congestion on one road does not affect that on another. In the conference version of this paper \cite{lazar2018poa}, we bound the price of anarchy and bicriteria for some nonseparable affine cost functions.
\end{remark}

The class of cost functions we consider are not monotone, meaning they do not necessarily satisfy \eqref{eq:monotonicity}, but are elementwise monotone, defined below:
\begin{definition}
	A cost function $c:R^{2n}_{\geq 0}\to R^{2n}_{\geq 0}$ is \emph{elementwise monotone} if it is nondecreasing in each of its arguments, \emph{i.e.} $\frac{\der c_i(z)}{\der z_j}\ge 0 \; \forall i,j \in {[}2n{]}$.
\end{definition}

%% file: bounds.tex
In this section we bound the price of anarchy and bicriteria of traffic networks under mixed autonomy. As established in the introduction, the PoA is in general unbounded in traffic networks with mixed autonomy. However, we can establish a bound for the PoA by parameterizing it using the parameters defined below.

\begin{definition}\label{def:asymmetry}
	The \emph{degree of asymmetry} on a road is the maximum ratio of road space utilized by a car of one type to a car of another type on the same road, while traveling at nominal velocity. The \emph{maximum degree of asymmetry}, $k$, is the maximum of the above quantity over all roads in the network. Formally, $k:=\max_{i\in \nidcs}\max(\lhdwyi/\shdwyi,\shdwyi/\lhdwyi)$.
\end{definition}

Note that we do not assume that one vehicle type affects congestion more than another type on all roads. For example, autonomous vehicles may require shorter headways than human-driven vehicles on highways but longer headways on neighborhood roads to maintain safety for pedestrians.

\begin{definition}
	The \emph{maximum polynomial degree}, denoted $\sigma$, for a road network with cost functions in the form \eqref{eq:cost_fn}, is the maximum degree of a polynomial denoting the cost on all roads in the network: $\sigma = \max_{i\in \nidcs}\sigma_i$.
\end{definition}

We use $\mathcal{C}_{k,\sigma}$ to denote the class of cost functions of the form \eqref{eq:cost_fn}, with maximum degree of asymmetry $k$ and maximum polynomial degree $\sigma$,  with cost functions using $m_i$ from either capacity mode 1 in \eqref{eq:capm1} or capacity model 2 in \eqref{eq:capm2}. Let
\begin{align}\label{def:xi}
	\xi(\sigma):=\sigma (\sigma+1)^{-\frac{\sigma+1}{\sigma}} \; .
\end{align}
Note that for $\sigma \ge 1$, $\xi(\sigma) < 1$. With this, we present our first bound.
\begin{theorem}\label{thm:poa}
	Consider a class of nonatomic congestion games with cost functions drawn from $\mathcal{C}_{k,\sigma}$. Let $z^{EQ}$ be an equilibrium and $z^{OPT}$ be a social optimum for this game. Then, 
	\begin{align*}
		C(z^{EQ})\le \frac{k^\sigma}{1-\xi(\sigma)}C(z^{OPT}) \; .
	\end{align*}
\end{theorem}
\begin{proof}
	
	Given any road cost function $c$ (and social cost $C$) and equilibrium $\zeq$, we define an aggregate cost function $\cagg$ (and social cost $\Cagg$) with corresponding equilibrium flow $\feq$, both parameterized by $\zeq$. This allows us to combine human-driven and autonomous flow into one flow type in the aggregate function so we can bound the price of anarchy for the aggregate cost function. We then find the relationship between the optimal routing for the aggregate cost function to that of the original cost function. Formally, the steps of the proof are:
	\begin{align}
	C(\zeq) &= \Cagg(\feq) \label{eq:aggregate} \\
	&\le \frac{1}{1-\xi(\sigma)}\Cagg(\fopt) \label{eq:bound_agg_poa} \\
	&\le \frac{1}{1-\xi(\sigma)}k^\sigma C(\zopt) \label{eq:bound_agg_opt} \; .
	\end{align}
	We begin by introducing the tool with which we bound the PoA in \eqref{eq:bound_agg_poa}. We then define $\cagg$ and $\feq$ such that \eqref{eq:aggregate} holds and show that $\feq$ is an equilibrium for $\cagg$. We discuss the structure of the tool used to bound the PoA and provide an intuitive explanation of how the chosen structure of $\cagg$ leads to a tighter PoA bound than an alternative choice. We then provide lemmas corresponding to inequality \eqref{eq:bound_agg_poa}, which bounds the price of anarchy of this new cost function, and \eqref{eq:bound_agg_opt}, which relates the social cost of optimal routing under $\cagg$ to that of the original cost function, $c$. We defer proofs of the lemmas to the appendix.
	
	We first introduce a general tool that we use for our results by extending the framework established by Correa \emph{et. al.} \cite{correa2008geometric}, which relies on the Variational Inequality to bound the price of anarchy. We use the following parameters:
	\begin{align}
	\beta(c,q) &:= \max_{z\in \mathbb{R}_{\ge 0}^{2N}} \frac{\langle c(q) - c(z), z\rangle}{\langle c(q),q \rangle} \; , \nonumber \\
	\beta(\mathcal{C})&:= \sup_{c\in\mathcal{C},q\in \feasrout} \beta(c,q) \; , \label{eq:beta_C}
	\end{align}
	where 0/0=0 by definition, and $\mathcal{C}$ is the class of network cost functions being considered. Then,
	\begin{lemma}\label{lma:correa}
		Let $z^{EQ}$ be an equilibrium of a nonatomic congestion game with cost functions drawn from a class $\mathcal{C}$ of elementwise monotone cost functions.	
		\begin{enumerate}[(a)]
			\item If $z^{OPT}$ is a social optimum for this game and $\beta(\mathcal{C}) < 1$, then
			\begin{align*}
			C(z^{EQ})\le (1-\beta(\mathcal{C}))^{-1}C(z^{OPT}) \; .
			\end{align*}
			\item If $q^{OPT}$ is a social optimum for the same game with $1+\beta(\mathcal{C})$ times as much flow demand of each type, then
			\begin{align*}
			C(z^{EQ})\le C(q^{OPT}) \; .
			\end{align*}
		\end{enumerate}
	\end{lemma}
	
	The lemma and proof are nearly identical to that of Correa \emph{et. al.} \cite{correa2008geometric}, extended to encompass nonmonotone, yet elementwise monotone, cost functions.
	
	We now explain our choice of $\cagg$ and $\fopt$ that yields \eqref{eq:aggregate} then provide an intuitive explanation for this choice. Recall that we define 
	\begin{align*}
	\zeq = \begin{bmatrix}
	\xeq_1 & \yeq_1 & \xeq_2 & \yeq_2 & \ldots & \xeq_n & \yeq_n
	\end{bmatrix}^T \; .
	\end{align*}
	
	We define a new flow vector which aggregates the regular and autonomous flows: $\feq = \xeq+\yeq$, where $\zeq\in \mathbb{R}^{2N}_{\ge 0}$ and $\xeq,\yeq,\feq \in \mathbb{R}^N_{\ge 0}$. We define a new cost function $\cagg$ that is a mapping from flow vector (with one flow for each road) to road latencies, \emph{i.e.} $\cagg: \mathbb{R}^N_{\ge 0} \rightarrow \mathbb{R}^N_{\ge 0}$. We define $\cagg$ so that it has the same road costs with flow $\feq$ as $c$ does with flow $\zeq$. Note, however, that $c$ is a mapping from flows, with two flow types per road, to road latencies, again with each road represented twice ($c:\mathbb{R}^{2N}_{\ge 0} \rightarrow \mathbb{R}^{2N}_{\ge 0}$). However, $\cagg$ represents each road once.
	
	We formally define $\cagg$, which depends on the equilibrium flow being considered, $\zeq$. This cost function is defined for both capacity models. In general terms, $\cagg$ adds the ``costly" type of vehicle flow first, then adds the ``less costly" vehicle flow. If $\shdwy_i\le \lhdwyi$,
	\begin{align}
	&\cagg_{i,1}(f_i) := \label{eq:agg_cost1} \\
	&\begin{cases}
	\ffli(1+\costscli(\frac{\lhdwyi f_i}{d_i})^{\sigma_i}) & f_i \leq \xeq_i \\
	\ffli(1+\costscli(\frac{\shdwyi f_i + (\lhdwyi - \shdwyi)\xeq_i }{d_i})^{\sigma_i}) & f_i >  \xeq_i
	\end{cases} \nonumber \\
	&\cagg_{i,2}(f_i) := \label{eq:agg_cost2} \\
	&\begin{cases} 
	\ffli(1+\costscli(\frac{\lhdwyi f_i}{d_i})^{\sigma_i}) & f_i \leq \xeq_i \\
	\ffli(1+\costscli(\frac{\lhdwyi f_i^2 - (\lhdwyi-\shdwyi)(f_i-\xeq_i)^2}{d_if_i})^{\sigma_i}) & f_i >  \xeq_i
	\end{cases} \nonumber
	\end{align}
	
	If $\shdwyi>\lhdwyi$, then swap $\shdwyi$ and $\lhdwyi$ above, and replace $\xeq_i$ with $\yeq_i$. In all cases, $\cagg_i(\feq_i) = c_j(\xeq_j,\yeq_j)$, where $j \in {[}2n{]}$ and $i=\lceil j/2 \rceil \in {[}n{]}$. Since the road latencies under $\cagg(\feq)$ are the same as under $c(\zeq)$, $\feq$ is an equilibrium for $\cagg$.
	
	\begin{figure}
		\centering
		\includegraphics[width=\linewidth]{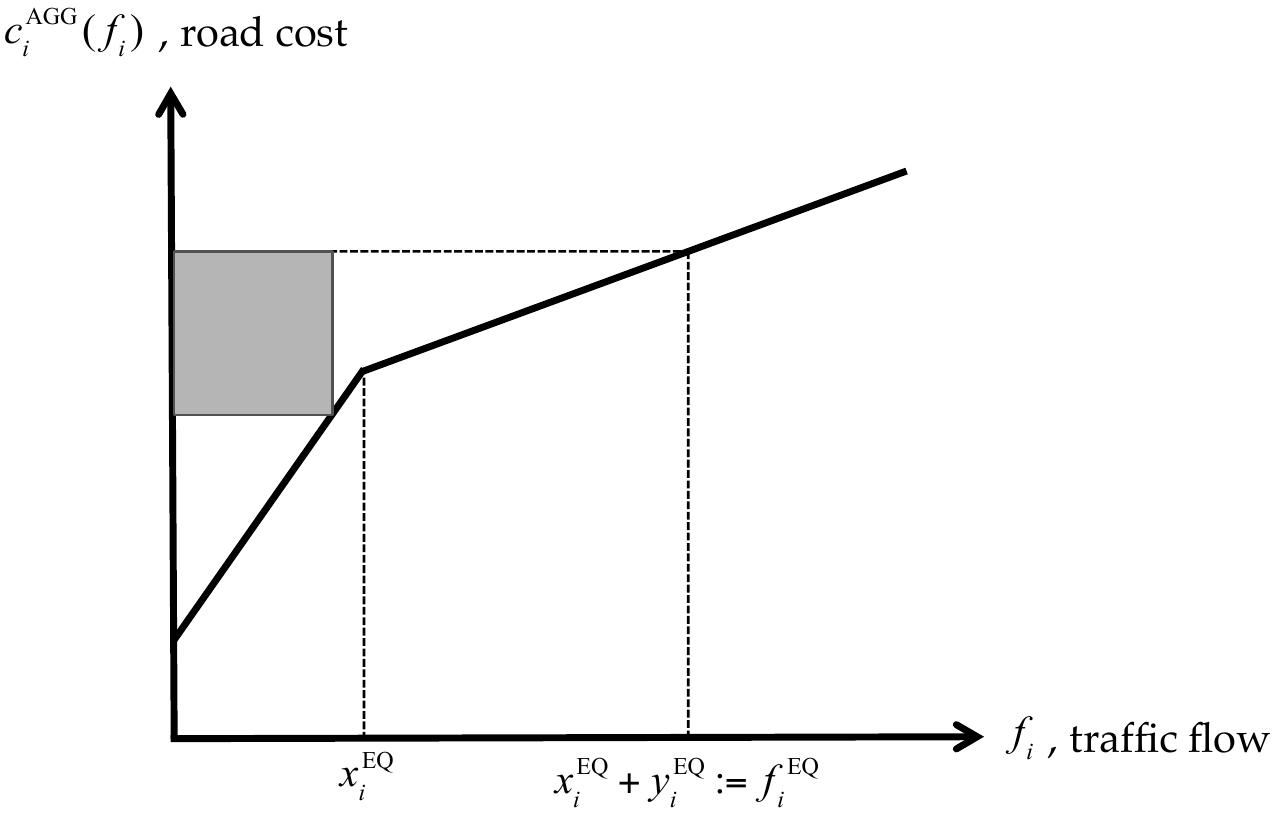}	
		\caption{Illustration of the geometric interpretation of the parameter $\beta(\mathcal{C}^\text{AGG})$ where $\mathcal{C}^\text{AGG}$ represents the class of aggregate cost functions. Parameter $\beta(\mathcal{C}^\text{AGG})$ is an upper bound on the ratio between the size of the shaded rectangle with the dashed rectangle. This is an upper bound over all choices of $c_i^{\text{AGG}}\in \mathcal{C}^\text{AGG}$ and $\xeq$ and $\yeq$ $\ge 0$.}
		\label{fig:agg_cost_fn}	
	\end{figure}
	
	To provide some intuition as to why we add the ``costly" vehicle type first, consider the affine case with the first capacity model. Correa \emph{et. al.} give a geometric interpretation of the parameter $\beta(\mathcal{C}^{\text{AGG}})$ when cost are separable, meaning road latency only depends on one element of the flow vector. They show that for any cost function drawn from $\mathcal{C}^{\text{AGG}}$, $\beta(\mathcal{C}^{\text{AGG}})$ provides an upper bound on the ratio of the area of a rectangle above the cost function curve to the area of a rectangle enclosing it, where the enclosing rectangle has one corner at the origin. See Figure~\ref{fig:agg_cost_fn} for an illustration.
	
	This interpretation provides the intuition that the more convex a function can be, the greater $\beta(\mathcal{C}^{\text{AGG}})$ can grow. Thus, to make our bound as tight as possible in our case, we add the costly vehicle type first. In the affine case with the first capacity model, this makes the class of cost functions concave. Then the element of this class that maximizes the size of the interior rectangle relative to the exterior rectangle minimizes the concavity of the function by setting $\xeq=0$ or $\yeq=0$. Thus, the PoA bound does not depend on the degree of asymmetry. Though this exact interpretation does not apply for $\sigma>1$ or for the second capacity model, the intuition is nonetheless useful.
	
	With this intuition, we now present inequalities \eqref{eq:bound_agg_poa} and \eqref{eq:bound_agg_opt} as a lemmas.
	\begin{lemma}\label{lma:agg_poa}
		Consider a nonatomic congestion game with road cost functions of the form \eqref{eq:agg_cost1} or \eqref{eq:agg_cost2}, with maximum polynomial degree $\sigma$. Then,
		\begin{align*}
		\Cagg(\feq) \le \frac{1}{1-\xi(\sigma)}\Cagg(\fopt)
		\end{align*}
		where $\xi(\sigma) = \sigma(\sigma+1)^{-\frac{\sigma+1}{\sigma}}$.
	\end{lemma}	
	
	\begin{lemma}\label{lma:agg_opt}
		Let $c$ be a cost function composed of road costs of the form \eqref{eq:cost_fn} with maximum degree of asymmetry $k$ and maximum polynomial degree $\sigma$. Let $\cagg$ be an aggregate cost function of $c$, as defined in \eqref{eq:agg_cost1} and \eqref{eq:agg_cost2}.	Let the flow vector $\zopt$ be a minimizer of $C$ and $\fopt$ be a minimizer of $\Cagg$, with $\sum_{i\in {[}2n{]}}\zopt_i = \sum_{i\in {[}n{]}}\fopt_i$. Then,
		\begin{align*}
		\Cagg(\fopt) \le k^\sigma C(\zopt) \; .
		\end{align*}
	\end{lemma}

	We prove Lemma \ref{lma:agg_poa} by bounding $\beta(\mathcal{C})$ for the class of aggregate cost functions and applying Lemma \ref{lma:correa}. We analyze the structures of $c$ and $\cagg$ to prove Lemma \ref{lma:agg_opt}. With these lemmas, the theorem is proved.
	
\end{proof}

Note that for $k=1$ (\emph{i.e.} no asymmetry), the price of anarchy bound simplifies to those in \cite{roughgarden2002selfishthesis, roughgarden2003topology, roughgarden2005selfish}. If the cost functions are affine and there is no asymmetry, this reduces to the classic $\frac{4}{3}$ bound \cite{roughgarden2002bad}. We characterize the tightness of this bound in the following corollary:
\begin{corollary}\label{cor:tightness}
	Given a maximum polynomial degree $\sigma$, the PoA bound is order-optimal with respect to the maximum degree of asymmetry $k$.
\end{corollary}

We provide an example proving the corollary in Section~\ref{sct:tightness}. When considering road networks with low asymmetry, we can establish another bound.
\begin{theorem}\label{thm:low_asym}
	Consider a class of nonatomic congestion games with cost functions drawn from $\mathcal{C}_{k,\sigma}$. Let $z^{EQ}$ be an equilibrium and $z^{OPT}$ a social optimum for this game. If $k\xi(\sigma)<1$, then 
	\begin{align*}
		C(z^{EQ})\le \frac{1}{1-k\xi(\sigma)}C(z^{OPT}) \; .
	\end{align*}
\end{theorem}

\begin{proof}
	To prove this theorem, instead of going through an aggregate cost function we directly find $\beta(\mathcal{C})$ for our class of cost functions and apply Lemma \ref{lma:correa}. We do this in two two lemmas: we first find a relationship between the parameter $\beta(c,v)$ and the road capacity model $m_i(x_i,y_i)$, then we bound the resulting expression.
	\begin{lemma}\label{lma:beta_bound}
		For cost functions of the form \eqref{eq:cost_fn}, the parameter $\beta(\mathcal{C})$ is bounded by
		\begin{align*}
		\beta(\mathcal{C}) \le \max_{i \in \nidcs,q_i,z_i\in \mathbb{R}_{\ge 0}^{2}}\frac{x_i + y_i}{v_i + w_i}(1 - (\frac{m_i(v_i,w_i)(x_i + y_i)}{m_i(x_i,y_i)(v_i + w_i)})^\sigma)
		\end{align*}
	\end{lemma}

	\begin{lemma}\label{lma:bound_beta_capms}
		For capacities of the forms \eqref{eq:capm1} or \eqref{eq:capm2}, 
		\begin{align*}
			\max_{i \in \nidcs,q_i,z_i\in \mathbb{R}_{\ge 0}^{2}}\frac{x_i + y_i}{v_i + w_i}(1 - (\frac{m_i(v_i,w_i)(x_i + y_i)}{m_i(x_i,y_i)(v_i + w_i)})^\sigma) \le k \xi(\sigma)
		\end{align*}
	\end{lemma}

	These lemmas, together with Lemma \ref{lma:correa}, prove the theorem as well as Theorem \ref{thm:bicriteria} below.

\end{proof}

Note that this bound may not necessarily be tighter in all regimes so our new PoA bound is $\min(\frac{k^\sigma}{1-\xi(\sigma)}, \frac{1}{1-k\xi(\sigma)})$. Though we cannot in closed form determine the region for which it is tighter, we can do so numerically. For example, for affine cost functions with $k=2$, $\frac{k^\sigma}{1-\xi(\sigma)}=\frac{8}{3}$ and $\frac{1}{1-k\xi(\sigma)}=2$. In Section~\ref{sct:tightness} we show via example that the bound in Theorem~\ref{thm:low_asym} is tight in this case.

The method used for establishing Theorem~\ref{thm:low_asym} also gives a bound on the bicriteria, stated in the following theorem.
\begin{theorem}\label{thm:bicriteria}
	Consider a class of nonatomic congestion games with cost functions drawn from $\mathcal{C}_{k,\sigma}$. Let $\zeq$ be an equilibrium for this game. If $\qopt$ is a social optimum for the same game with $1+ k \xi(\sigma)$ times as much flow demand of each type, then 
	\begin{align*}
		C(\zeq)\le C(\qopt) \; .
	\end{align*}
\end{theorem}

To give an example, if road delays are described by polynomials of degree 4, and the maximum asymmetry between the spacing of platooned and nonplatooned vehicles is 3, then the the cost of selfishly routing vehicles will be less than optimally routing $1+3 \, \xi(4) \approx 2.61$ times as much vehicle flow of each type.

%% file: discussion_examples.tex
\begin{figure}
	\centering
	\begin{subfigure}[b]{1\linewidth}
		\centering
			\begin{center}
				\begin{tikzpicture}[->, >=stealth', auto, semithick, node distance=5cm]
				\tikzstyle{every state}=[fill=white,draw=black,thick,text=black,scale=1]
				\node[state]    (0)               {$s$};
				\node[state]    (1)[right of=0]   {$t$};
				\path
				(0) edge[bend left]		node[above]{$c_1(x,y) = (kx+y)^\sigma$}     (1)
				(0) edge[bend right]		node[below]{$c_2(x,y) = (x+ky)^\sigma$}     (1);
				\end{tikzpicture}
			\end{center}
	\end{subfigure}
	\begin{subfigure}[b]{1\linewidth}
		\centering
		\includegraphics[width=\linewidth]{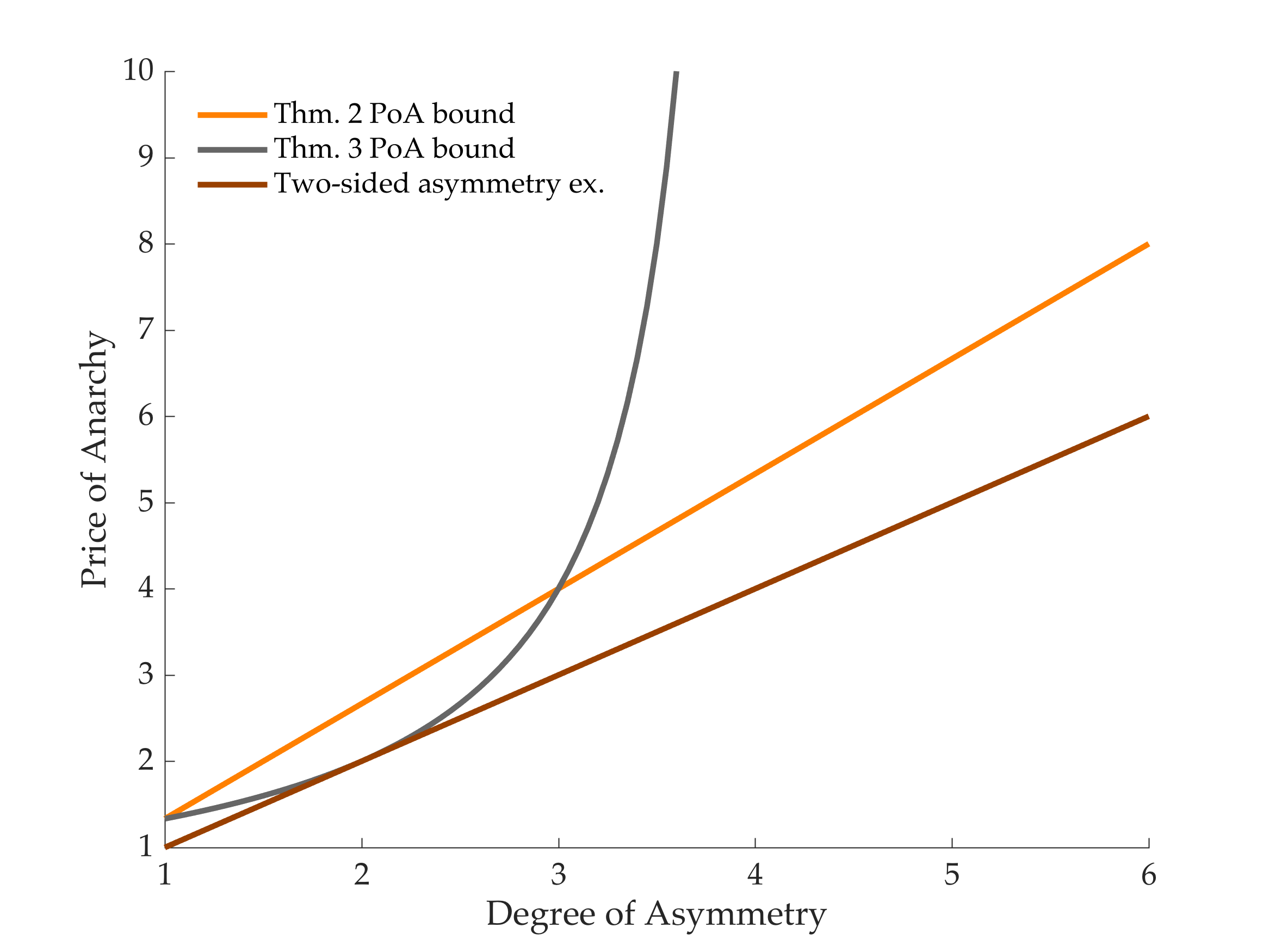}	
	\end{subfigure}
	\caption{Example road network with two sided asymmetry (above) and comparison of price of anarchy with upper bound (below). One unit of human-driven and one unit of autonomous flow cross from node $s$ to node $t$.}	
	\label{fig:results_two_sided}
\end{figure}

\begin{figure}
	\centering
	\begin{subfigure}[b]{1\linewidth}
		\centering
		\begin{center}
			\begin{tikzpicture}[->, >=stealth', auto, semithick, node distance=5cm]
			\tikzstyle{every state}=[fill=white,draw=black,thick,text=black,scale=1]
			\node[state]    (0)               {$s$};
			\node[state]    (1)[right of=0]   {$t$};
			\path
			(0) edge[bend left]		node[above]{$c_1(x,y) = 1$}     (1)
			(0) edge[bend right]	node[below]{$c_2(x,y) = \frac{k}{\sqrt{k}+1}x+\frac{1}{\sqrt{k} + 1}y$}     (1);
			\end{tikzpicture}
		\end{center}
	\end{subfigure}
	\begin{subfigure}[b]{1\linewidth}
		\centering
		\includegraphics[width=\linewidth]{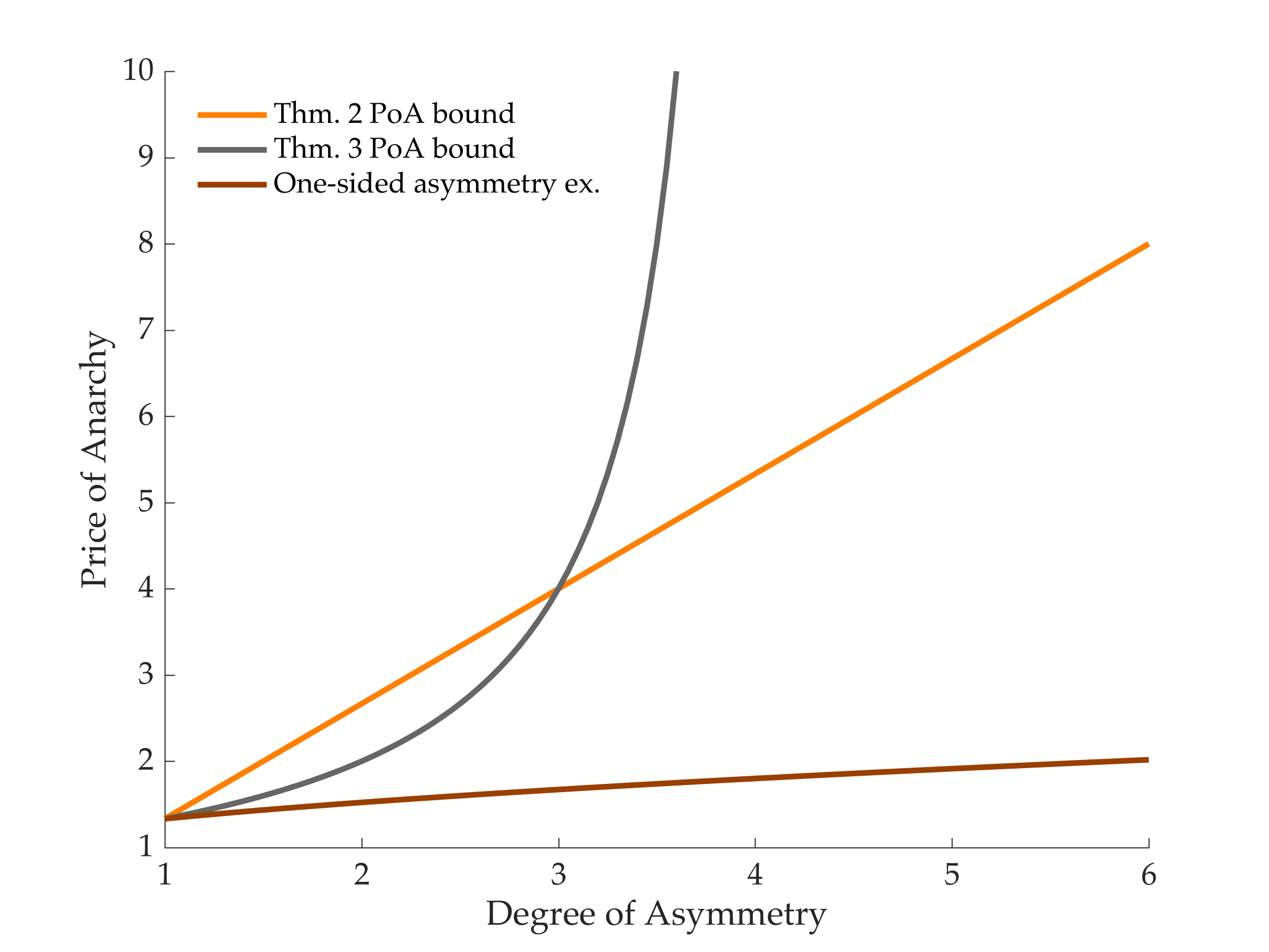}
	\end{subfigure}
	\caption{Example road network with one sided asymmetry (above) and comparison of price of anarchy with upper bound for $\sigma=1$ (below). $\frac{1}{\sqrt{k}}$ units of human-driven flow and one unit of autonomous flow cross from node $s$ to $t$.}	
	\label{fig:results_one_sided}
\end{figure}
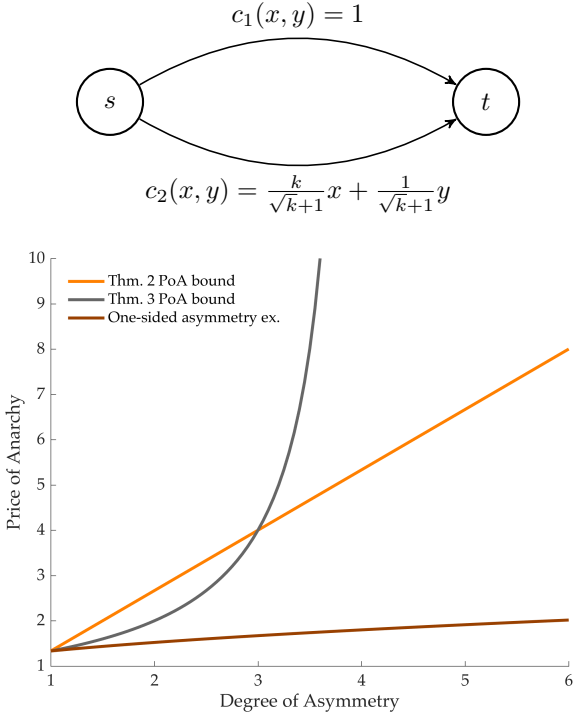

\begin{figure}
	\centering
	\includegraphics[width=\linewidth]{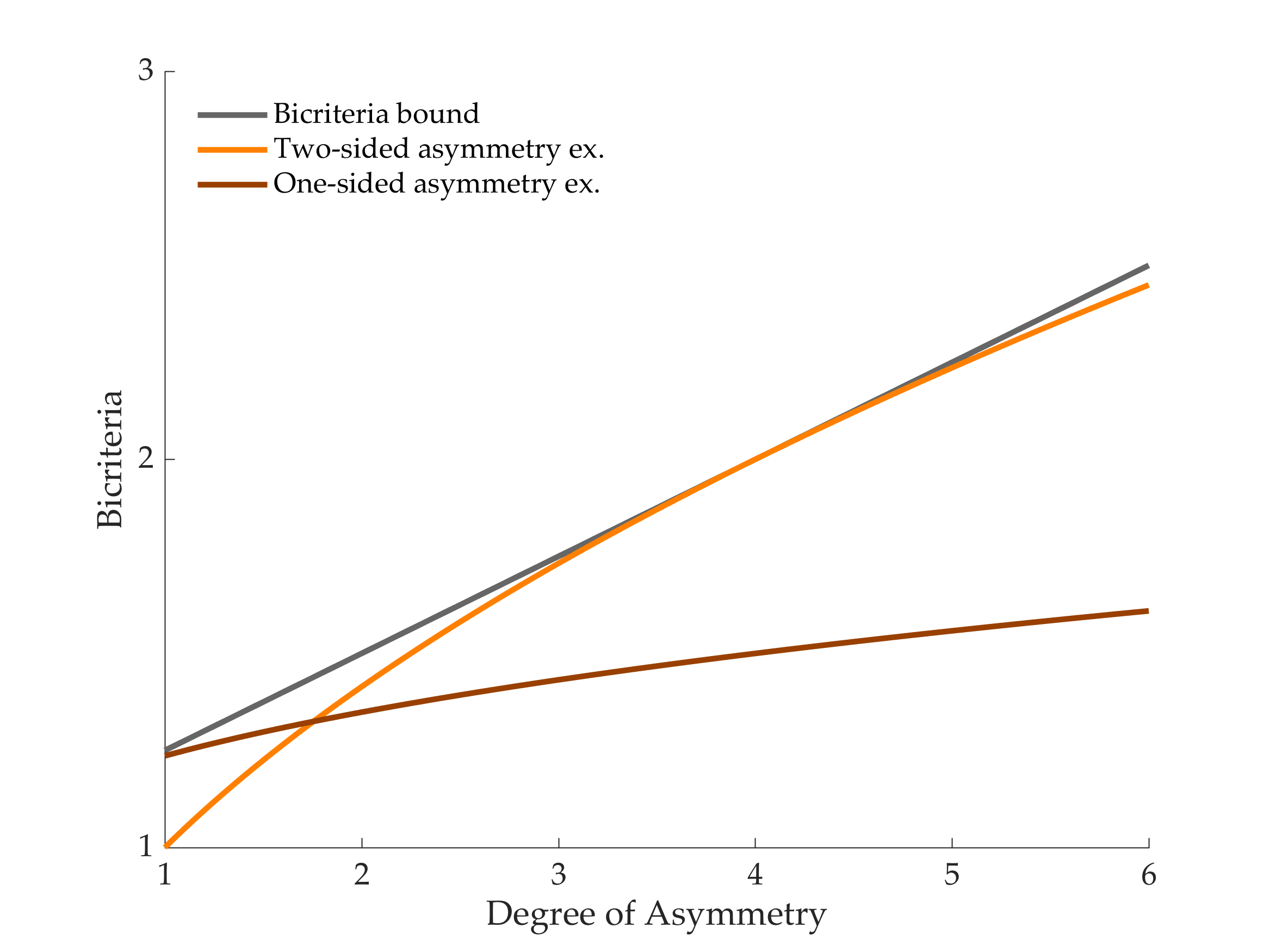}	
	\caption{Bicriteria of networks in Figures \ref{fig:results_two_sided} (with $\sigma=1$) and \ref{fig:results_one_sided}, compared with the bicriteria bound.}
	\label{fig:bicriteria}	
\end{figure}

In this section we provide examples that give a lower bound on the price of anarchy for this class of networks and serve to illustrate the tightness of the bounds. We discuss notions of one-sided and two-sided asymmetry: a network has one sided asymmetry if $\shdwyi \le \lhdwyi~\forall i \in \nidcs$ (human-driven cars always contribute more to congestion than autonomous cars) or $\shdwyi \ge \lhdwyi~\forall i \in \nidcs$ (human-driven cars always contribute less to congestion than autonomous cars); otherwise the network has two-sided asymmetry. We provide two example networks, one with two-sided asymmetry (Fig.~\ref{fig:results_two_sided}) and one with one-sided asymmetry (Fig.~\ref{fig:results_one_sided}). We compare the price of anarchy and bicriteria in those networks to the upper bounds established earlier.

\begin{example}\label{ex:two_sided_network}
	Consider the network of parallel roads in Fig.~\ref{fig:results_two_sided}, where one unit of human-driven and one unit of autonomous flow wish to cross from node $s$ to $t$. The roads have costs $c_1(x,y) = (kx+y)^\sigma$ and $c_2(x,y) = (x+ky)^\sigma$. In worst-case equilibrium, all human-driven cars are on the top road and all autonomous cars are on the bottom road. In the best case, these routing are reversed. This yields a PoA of $k^\sigma$. To find the bicriteria, we calculate how much traffic, optimally routed, yields a cost equal to $2k^\sigma$, the cost of the original traffic volume at worst-case equilibrium. We find that $k^{\frac{\sigma}{\sigma + 1}}$ as much traffic, optimally routed, yields this same cost.
\end{example}

\begin{example}\label{ex:one_sided_network}
	Consider the network of parallel roads in Figure \ref{fig:results_one_sided}, where $\frac{1}{\sqrt{k}}$ unit of human-driven and one unit of autonomous flow wish to cross from node $s$ to $t$. The roads have costs $c_1(x,y) = 1$ and $c_2(x,y) = \frac{k}{\sqrt{k}+1}x+\frac{1}{\sqrt{k} + 1}y$. At equilibrium, all vehicles take the bottom road; optimally routed, human-driven vehicles takes the top road and autonomous vehicles take  the bottom. This yields a price of anarchy of $1+\frac{k}{2\sqrt{k}+1}$. Calculations similar to that in Example \ref{ex:two_sided_network} yield a bicriteria of $\frac{(-1 + \sqrt{1 + 4\sqrt{k}})(1+\sqrt{k})}{2\sqrt{k}}$.
\end{example}

For affine cost functions, $\sigma=1$ so $\xi=1/4$. The PoA bound is then $\min(\frac{4}{4-k},\frac{4}{3}k)$ and the bicriteria bound is $1+k/4$. With $\sigma=1$, the first example has PoA $k$ and bicriteria $\sqrt{k}$, and the second example has PoA of order $\sqrt{k}$ and bicriteria of order $k^{1/4}$. Accordingly, the first example shows that the PoA bound is tight for $k=2$ and the bicriteria bound is tight for $k=4$. Note that a realistic range for $k$ is between $1$ and $4$.

Further, for affine cost functions, the bound in Theorem~\ref{thm:low_asym} is tighter than that of Theorem~\ref{thm:poa} when the degree of asymmetry is low. However, the bound in Theorem~\ref{thm:poa} scales much better for high degrees of asymmetry. This effect is accentuated for cost functions that have higher order polynomials -- the regime for which the bound in Theorem~\ref{thm:low_asym} is tighter shrinks as the maximum polynomial degree grows.

As stated in Corollary \ref{cor:tightness}, our bound is order-optimal with respect to the maximum degree of asymmetry, $k$.  Comparing the bound ($\frac{k^\sigma}{1-\xi(\sigma)}$) with the PoA in Example \ref{ex:two_sided_network} ($k^\sigma$) shows that for a fixed $\sigma$, the PoA upper bound is within a constant factor of the lower bound, implying that the upper bound is order-optimal in $k$.

It is also worth noting that under the construction used in Theorems \ref{thm:low_asym}, the bicriteria is related to the price of anarchy through the quantity $\beta(\mathcal{C})$ \cite{correa2008geometric}. Observe that Example \ref{ex:two_sided_network} provides a bicriteria of 2 for $k=4$, implying $\beta(\mathcal{C}_4)\ge 1$. Since the PoA is greater than or equal to $\frac{1}{1-\beta(\mathcal{C})}$, this mechanism cannot bound the PoA for $k\ge 4$. This leads us to rely on the mechanism developed for Theorem~\ref{thm:poa} for networks with large asymmetry.

%% file: conclusion.tex
In this paper we present a framework, similar to a congestion game, for considering traffic networks with mixed autonomy. To do so we present two models for the capacity of roads with mixed autonomy, each corresponding to an assumption about the technological capabilities of autonomous vehicles, and we define a class of road latency functions that incorporates these capacity models. Using this framework, we develop two methods of bounding the price of anarchy and show that these bounds depend on the degree of the polynomial describing latency and the difference in the degree to which platooned and nonplatooned vehicles occupy space on a road. In addition we present a bound on the bicriteria, another measure of inefficiency due to selfish routing. We present examples showing these bounds are tight in some cases and recover classical bounds when human-driven and autonomous vehicles affect congestion the same way. Moreover, we show that our PoA bound is order-optimal with respect to the degree to which vehicle types differently affect congestion.

Some limitations of the work are as follows. The capacity models presented assume that vehicle types are determined as a result of a Bernoulli process; a more general capacity model could incorporate autonomous vehicles that actively rearrange themselves as to form platoons. Further, the road latency model assumes that latency increases with increasing flow. This does not incorporate the notion of congestion described by the fundamental diagram of traffic, in which a road can have both low flow and high latency if vehicle density is high. Similarly, autonomous vehicles can affect congestion in ways not limited to platooning. In addition, our proposed latency function considers only the effect of a vehicle on the road upon which it travels; a more general latency function would consider interaction between roads. Finally, the PoA bound is not shown to be tight but is order-optimal in the degree of asymmetry $k$, and a future work could aim to close this gap. Nonetheless, this is a framework that can be used in the future for studying traffic networks in mixed autonomy.

%% file: proofs.tex
\subsection{Proof of Lemma \ref{lma:correa}}\label{pf:correa}

To prove part (a),
\begin{align}
\langle c(\qeq),z \rangle &= \langle c(z), z\rangle + \langle c(\qeq)-c(z), z\rangle \nonumber \\
&\le \langle c(z), z \rangle + \beta(c,\qeq)\langle c(\qeq),\qeq \rangle \nonumber \\
&\le C(z) + \beta(\mathcal{C})C(\qeq) \label{eq:lemma41}
\end{align}	
and by the Variational Inequality, $C(\qeq) \le \langle c(\qeq),z \rangle$ for any feasible routing $z$. Completing the proof requires that $\beta(\mathcal{C}) \le 1$, then replace the generic $z$ with $\zopt$.

To prove part (b), element-wise monotonicity implies the feasibility of $(1+\beta(\mathcal{C}))^{-1}\qopt$, which routes the same volume of traffic as $\zeq$. Using \eqref{eq:VI},
\begin{align}
\langle c(\zeq),\zeq \rangle \le \langle c(\zeq),(1+\beta(\mathcal{C}))^{-1}\qopt \rangle \; . \label{eq:thm1ln0}
\end{align}
Then,
\begin{align}
C(\zeq) &= (1+\beta(\mathcal{C}))\langle c(\zeq),\zeq \rangle \nonumber \\
& \quad - \beta(\mathcal{C})\langle c(\zeq), \zeq \rangle \label{eq:thm1ln1} \\
&\le (1+\beta(\mathcal{C}))\langle c(\zeq),(1+\beta(\mathcal{C}))^{-1}\qopt \rangle \nonumber \\
& \quad - \beta(\mathcal{C})\langle c(\zeq), \zeq \rangle \label{eq:thm1ln2}\\
&\le C(\qopt) \; , \label{eq:thm1ln3}
\end{align}
where \eqref{eq:thm1ln2} uses \eqref{eq:thm1ln0} and \eqref{eq:thm1ln3} uses \eqref{eq:lemma41}.

\subsection{Proof of Lemma \ref{lma:agg_poa}}\label{pf:agg_poa}

We use Lemma~\ref{lma:correa} and bound $\beta(\mathcal{C}^\text{AGG}_{\sigma})$, where $\mathcal{C}^\text{AGG}_{\sigma}$ denotes the set of aggregate cost functions with maximum polynomial degree $\sigma$. First we will show that 
\begin{align*}
\beta(\mathcal{C}^\text{AGG}_{\sigma}) \le \max_{i \in \nidcs, f_i, g_i\ge 0} \frac{f_i}{g_i}(1-(\frac{f_i}{g_i})^{\sigma_i}) \; .
\end{align*}

For the remainder of the proof of the lemma, we drop the aggregate superscript from aggregate cost functions. Then,	
\begin{align*}
\beta(\mathcal{C}^{\text{AGG}}_{\sigma}) &:= \sup_{c\in \mathcal{C}^\text{AGG}_{\sigma}, g \in \feasrout} \max_{f \in \mathbb{R}^n_{\ge 0}}\frac{\langle c(g)-c(f),f \rangle}{\langle c(g),c(g) \rangle} \\
&\le \sup_{c_i\in \mathcal{C}^\text{AGG}_{\sigma}, g_i \ge 0} \max_{f_i \ge 0}\frac{f_i}{g_i}(1-\frac{c_i(f_i)}{c_i(g_i)}) \; .
\end{align*}
We will show that $\beta(\mathcal{C}^\text{AGG}_{\sigma})\le \max_{i \in \nidcs, f_i, g_i\ge 0}\frac{f_i}{g_i}(1-(\frac{f_i}{g_i})^{\sigma_i})$ by showing that $\frac{c_i(f_i)}{c_i(g_i)}\ge (\frac{f_i}{g_i})^{\sigma_i}$ for either capacity model. We will show this for a road in which $\shdwy_i\le \lhdwyi$, though with the alterations discussed above, the same can be done for a road on which $\shdwy_i > \lhdwyi$.

Note that we are considering the supremum of $g$ in the set of feasible routings for $c$. Accordingly, we need to consider cases in which $g_i\le \xeq_i$, as well as $g_i > \xeq_i$. However, $g_i$ must be greater than $f_i$, as $\beta(\mathcal{C}^{\text{AGG}}_{k,\sigma})\ge 0$ and $c_i$ is nondecreasing. We therefore do not need to consider cases in which $g_i\le \xeq_i$ and $f_i>\xeq_i$.

All together, we must bound six cases:
\begin{enumerate}
	\item capacity model 1, $g_i \le \xeq_i$ and $f_i \le \xeq_i$,
	\item capacity model 1, $g_i > \xeq_i$ and $f_i \le \xeq_i$,
	\item capacity model 1, $g_i > \xeq_i$ and $f_i > \xeq_i$,
	\item capacity model 2, $g_i \le \xeq_i$ and $f_i \le \xeq_i$,
	\item capacity model 2, $g_i > \xeq_i$ and $f_i \le \xeq_i$,
	\item capacity model 2, $g_i > \xeq_i$ and $f_i > \xeq_i$,
\end{enumerate}
We will show that in all cases, $\frac{c_i(f_i)}{c_i(g_i)}\ge (\frac{f_i}{g_i})^{\sigma_i}$.

In the first case,
\begin{align}
\frac{c_{i,1}(f_i)}{c_{i,1}(g_i)} &= \frac{1+\costscli (\frac{\lhdwyi f_i}{d_i})^{\sigma_i}}{1+\costscli (\frac{\lhdwyi g_i}{d_i})^{\sigma_i}} \nonumber \\
&\ge (\frac{f_i}{g_i})^{\sigma_i} \label{eq:simp1} \; ,
\end{align}
where \eqref{eq:simp1} follows from $\frac{c_{i,1}(f_i)}{c_{i,1}(g_i)}\le 1$. In the following cases we perform the same operation without comment.

In the second case,
\begin{align}
\frac{c_{i,1}(f_i)}{c_{i,1}(g_i)} &\ge (\frac{h_i f_i}{\shdwyi g_i + (\lhdwyi - \shdwyi)\xeq_i})^{\sigma_i} \nonumber \\
& \ge (\frac{h_i f_i}{\shdwyi g_i})^{\sigma_i} \nonumber \\
& \ge (\frac{f_i}{g_i})^{\sigma_i} \label{eq:simp2} \; ,
\end{align}
where \eqref{eq:simp2} follows from $\shdwyi\le \lhdwyi$.

In the third case,
\begin{align}
\frac{c_{i,1}(f_i)}{c_{i,1}(g_i)} &\ge (\frac{\shdwyi f_i + (\lhdwyi - \shdwyi)\xeq_i}{\shdwyi g_i + (\lhdwyi - \shdwyi)\xeq_i})^{\sigma_i} \nonumber \\
&\ge (\frac{f_i}{g_i})^{\sigma_i} \label{eq:simp3} \; ,
\end{align}
where \eqref{eq:simp3} follows from $\lhdwyi \ge \shdwyi$ and $f_i\le g_i$.

The fourth case is equivalent to the first case. In the fifth case,
\begin{align}
\frac{c_{i,2}(f_i)}{c_{i,2}(g_i)} &\ge (\frac{\lhdwyi f_i g_i}{ \lhdwyi (g_i)^2 - (\lhdwyi-\shdwyi)(g_i-\xeq_i)^2 } )^{\sigma_i} \nonumber \\
&\ge (\frac{\lhdwyi f_i g_i}{ \lhdwyi (g_i)^2 } )^{\sigma_i} = (\frac{ f_i}{ g_i } )^{\sigma_i} \label{eq:simp4} \; ,
\end{align}
where \eqref{eq:simp4} results from $\lhdwyi \ge \shdwyi$.

In the sixth case,
\begin{align}
\frac{c_{i,2}(f_i)}{c_{i,2}(g_i)} &\ge (\frac{\lhdwyi (f_i)^2 - (\lhdwyi-\shdwyi)(f_i-\xeq_i)^2}{ \lhdwyi (g_i)^2 - (\lhdwyi-\shdwyi)(g_i-\xeq_i)^2 }\times \frac{g_i}{f_i})^{\sigma_i} \nonumber \\
&\ge (\frac{\lhdwyi (f_i)^2}{ \lhdwyi (g_i)^2}\times \frac{g_i}{f_i})^{\sigma_i} = (\frac{f_i}{ g_i})^{\sigma_i} \label{eq:simp5} \; ,
\end{align}
where \eqref{eq:simp5} results from $(\frac{f_i}{g_i})^2 \ge (\frac{f_i-\xeq_i}{g_i-\xeq_i})^2$, as $g_i \ge f_i$.

Now that we have shown that $\frac{c_{i}(f_i)}{c_{i}(g_i)}\ge (\frac{f_i}{g_i})^{\sigma_i}$ in all cases, we find that
\begin{align*}
\beta(\mathcal{C}^{\text{AGG}}_{\sigma}) &\le \max_{i\in \nidcs,f_i,g_i\ge 0} \frac{f_i}{g_i}- (\frac{f_i}{g_i})^{\sigma_i+1} \; .
\end{align*}

As this expression is concave with respect to $f_i$, to maximize this with respect to $f_i$, we set the derivative of this expression with respect to $f_i$ to 0:
\begin{align*}
&\frac{1}{g_i} - (\sigma + 1)\frac{(f^*_i)^{\sigma}}{(g_i)^{\sigma+1}} = 0 \\
&\implies f_i^*=(g_i)(\sigma+1)^{-\frac{1}{\sigma}} \; .
\end{align*}

Plugging this in, 
\begin{align*}
\beta(\mathcal{C}^{\text{AGG}}_{\sigma}) &\le (\sigma+1)^{-\frac{1}{\sigma}}(1-\frac{1}{\sigma+1})  \\
&= \sigma(\sigma+1)^{-\frac{\sigma+1}{\sigma}}
\end{align*}

This, combined with Lemma \ref{lma:correa} and the definition of $\xi(\sigma)$, completes the proof of the lemma.

\subsection{Proof of Lemma \ref{lma:agg_opt}}\label{pf:agg_opt}

As before, let  
\begin{align*}
z = \begin{bmatrix} x_1 & y_1 & x_2& y_2 & \ldots & x_n & y_n \end{bmatrix}^T \; .
\end{align*} Then,
\begin{align*}
k^\sigma &c_i(\xopt_i,\yopt_i) \\
&\ge \max(c_i(\xopt_i+\yopt_i,0), c_i(0, \xopt_i+\yopt_i)) \\
&\ge \cagg_i(\xopt_i+\yopt_i) \; .
\end{align*}
and by definition of $\fopt$, $\Cagg(f)\ge \Cagg(\fopt)$ for any feasible vector $f$ with $\sum_{i \in \nidcs}\fopt_i = \sum_{i \in \nidcs}f_i$, so
\begin{align*}
k^\sigma C(\zopt) &\ge \sum_{i \in \nidcs}(\xopt_i+\yopt_i)\cagg_i(\xopt+\yopt) \\
&\ge \Cagg(\fopt) \; .
\end{align*}

\subsection{Proof of Lemma \ref{lma:beta_bound}}\label{pf:beta_bound}

Using \eqref{eq:beta_C},
\begin{align}
&\beta(c,q) \nonumber \\
&= \max_{z_i \in \mathbb{R}_{\ge 0}^{2n}} \frac{\sum_{i\in \nidcs}\ffli \costscli {[} (\frac{v_i+w_i}{m_i(v_i,w_i)})^{\sigma_i}-(\frac{x_i+y_i}{m_i(x_i,y_i)})^{\sigma_i}{]}(x_i+y_i)}{\sum_{i\in \nidcs}\ffli {[} 1 + \costscli(\frac{v_i+w_i}{m_i(v_i,w_i)})^{\sigma_i} {]}(v_i+w_i)} \nonumber \\
&\le \max_{i\in \nidcs, z_i\in \mathbb{R}_{\ge 0}^{2}} \frac{\costscli{[}(\frac{v_i+w_i}{m_i(v_i,w_i)})^{\sigma_i}-(\frac{x_i+y_i}{m_i(x_i,y_i)})^{\sigma_i}{]}(x_i+y_i)}{{[}1+\costscli(\frac{v_i+w_i}{m_i(v_i,w_i)})^{\sigma_i}{]}(v_i+w_i)} \nonumber \\
&\le \max_{i\in \nidcs, z_i\in \mathbb{R}_{\ge 0}^{2}}\frac{[(\frac{v_i+w_i}{m_i(v_i,w_i)})^{\sigma_i}-(\frac{x_i+y_i}{m_i(x_i,y_i)})^{\sigma_i}](x_i+y_i)}{(\frac{v_i+w_i}{m_i(v_i,w_i)})^{\sigma_i}(v_i+w_i)} \nonumber \\
&= \max_{i \in \nidcs, z_i\in \mathbb{R}_{\ge 0}^{2}}\frac{x_i + y_i}{v_i + w_i}(1 - (\frac{m_i(v_i,w_i)(x_i + y_i)}{m_i(x_i,y_i)(v_i + w_i)})^{\sigma_i}) \nonumber \\
&\le \max_{i \in \nidcs, z_i\in \mathbb{R}_{\ge 0}^{2}}\frac{x_i + y_i}{v_i + w_i}(1 - (\frac{m_i(v_i,w_i)(x_i + y_i)}{m_i(x_i,y_i)(v_i + w_i)})^\sigma) \; , \label{eq:beta_simpl}
\end{align}
where $\beta(c,q)\ge 0$ implies \eqref{eq:beta_simpl}, allowing us to consider only the maximum allowable degree of polynomial.

\subsection{Proof of Lemma \ref{lma:bound_beta_capms}}\label{pf:bound_beta_capms}

For capacity model 1:
\begin{align}
&\beta(c,q) \\
&\le \max_{i \in \nidcs, z_i\in \mathbb{R}_{\ge 0}^{2}}\frac{x_i + y_i}{v_i + w_i}(1 - \nonumber \\
&\qquad \qquad (\frac{(h_i-(h_i-\bar{h}_i)(\frac{y_i}{x_i+y_i}))(x_i + y_i)}{(h_i-(h_i-\bar{h}_i)(\frac{w_i}{v_i+w_i}))(v_i + w_i)})^\sigma) \nonumber \\
&= \max_{i \in \nidcs, z_i\in \mathbb{R}_{\ge 0}^{2}}\frac{x_i + y_i}{v_i + w_i}(1 - (\frac{h_i x_i + \bar{h}_i y_i}{h_i v_i + \bar{h}_i w_i})^\sigma) \nonumber \\
&\le \max_{z_i\in \mathbb{R}_{\ge 0}^{2}}\frac{x_i + y_i}{v_i + w_i}(1 - (\frac{k x_i + y_i}{k v_i + w_i})^\sigma) \label{eq:beta_simpl_capm1} \; , \\
&:= f(x_i,y_i,v_i,w_i) \; . \nonumber
\end{align}

In \eqref{eq:beta_simpl_capm1} we use the Definition \ref{def:asymmetry} of the maximum degree of asymmetry. For ease of notation, we drop the subscripts for $f(x,y,v,w)$.

We now investigate this expression more closely, and show that the maximum of this expression with respect to $x$ and $y$ occurs at either $x=0$ or $y=0$. We do this by showing there are no critical points with $x>0$ and $y>0$, and that outside of a finite region, the function is decreasing with respect to both $x$ and $y$.

First we show that there exist no critical points, meaning points for which $\frac{\der f}{\der x} = \frac{\der f}{\der y} = 0$, for $k>1$. We have
\begin{align*}
\frac{\der f}{\der x} &= \frac{kx + y - (\frac{k x + y}{k v+w})^\sigma(k x + y + k \sigma (x+y))}{(v+w)(kx+y)} \\
\frac{\der f}{\der y} &= \frac{kx + y - (\frac{k x + y}{k v+w})^\sigma(k x + y + \sigma (x+y))}{(v+w)(kx+y)} \; .
\end{align*}
Since $k>1$, we conclude that $\frac{\der f}{\der x} \neq \frac{\der f}{\der y}$ for $x>0$ and $y>0$.

To show the second component, we see that
\begin{align*}
\frac{\der f}{\der x} & \le \frac{1}{v+w}(1-(\frac{kx+y}{kv+w})^\sigma) \quad \text{and} \\
\frac{\der f}{\der y} & \le \frac{1}{v+w}(1-(\frac{kx+y}{kv+w})^\sigma) \; .
\end{align*}
Therefore, for the region $y > kv+w-kx$, the function is decreasing with $x$ and $y$.

These two facts together imply the maximum of $f$ in the first quadrant lies on either the $x$ or $y$ axis. Checking these two candidate functions,
\begin{align*}
f(x,0,v,w) = \frac{x(1-(\frac{k x}{k v + w})^\sigma)}{v+w} \\
f(0,y,v,w) = \frac{y(1-(\frac{y}{k v + w})^\sigma)}{v + w}
\end{align*}
These functions are concave with respect to x and y, with minima at $x^* = (\frac{1}{\sigma+1})^{1/\sigma}\frac{k v + w}{k}$ and $y^* = (\frac{1}{\sigma+1})^{1/\sigma}(k v + w)$ respectively. When plugging these in, we find that the solution along the y-axis is greater, so
\begin{align*}
\max_{x,y\ge 0}f(x,y,v,w)= \sigma (\frac{1}{\sigma+1})^{1+1/ \sigma}\frac{k v + w}{v+w} \; .
\end{align*}
We then find  that when restricted to capacity model 1,
\begin{align*}
\beta(\mathcal{C}_{k,\sigma}) &\le \max_{v,w \le 0} \sigma (\frac{1}{\sigma+1})^{1+1/ \sigma}\frac{k v + w}{v+w} \\
&\le k \sigma (\frac{1}{\sigma+1})^{1+1/ \sigma} \\
&= k \xi(\sigma)  \; .
\end{align*}

For capacity model 2:
\begin{align*}
&\beta(c,q) \\
&\le \max_{z_i\in \mathbb{R}_{\ge 0}^{2}}\frac{x_i+y_i}{v_i+w_i}-(\frac{kx_i^2 + 2kx_iy_i + y_i^2}{kv_i^2 + 2kv_iw_i+ w_i^2})^\sigma(\frac{v_i+w_i}{x_i+y_i})^{\sigma-1} \\
& := g(x_i,y_i,v_i,w_i)
\end{align*}
We again find that $k>1 \implies \frac{\der g}{\der x} \neq \frac{\der g}{\der y}$. Further,
\begin{align*}
\frac{\der f}{\der x}&\le \frac{1}{v+w}(1-((\frac{v+w}{x+y})\frac{kx^2 + 2kxy + y^2}{kv^2 + 2kvw + w^2})^\sigma \quad \text{and} \\
\frac{\der f}{\der y}&\le \frac{1}{v+w}(1-((\frac{v+w}{x+y})\frac{kx^2 + 2kxy + y^2}{kv^2 + 2kvw + w^2})^\sigma \; ,
\end{align*}
so
\begin{align*}
&\frac{kx^2 + 2kxy + y^2}{x+y}<\frac{kv^2 + 2kvw + w^2}{v+w} \\
&\quad \implies \frac{\der f}{\der x}<0 \; \text{and} \; \frac{\der f}{\der y}<0 \; .
\end{align*}
Using the same reasoning as above, we now just search the x- and y-axes. As above, $g(x, 0, v, w)$ is concave with respect to $x$ and $g(0,y,v,w)$ is concave with respect to $y$, with maxima at $x=(\frac{1}{\sigma+1})^{1/\sigma}\frac{kv^2 + 2kvw + w^2}{k(v+w)}$ and $y=(\frac{1}{\sigma+1})^{1/\sigma}\frac{kv^2 + 2kvw + w^2}{v+w}$, respectively. Comparing these, we find the maximum is the latter, so
\begin{align*}
\max_{x,y\ge 0}g(x,y,v,w) = \sigma(\frac{1}{\sigma+1})^{1+1/\sigma}\frac{kv^2 + 2kvw + w^2}{(v+w)^2} \; .
\end{align*}
Then, for capacity model 2,
\begin{align*}
\beta(\mathcal{C}_{k,\sigma}) &\le \max_{v,w \le 0} \sigma (\frac{1}{\sigma+1})^{1+1/ \sigma}\frac{kv^2 + 2kvw + w^2}{(v+w)^2} \\
&\le \sigma (\frac{1}{\sigma+1})^{1+1/ \sigma}k \\
&= k \xi(\sigma) \; .
\end{align*}
Together, this shows that regardless of capacity model, $\beta(\mathcal{C}_{k,\sigma})\le k  \xi(\sigma)$. The application of Lemma \ref{lma:correa} completes the proof.

%% file: bios.tex
\begin{IEEEbiographynophoto}{Daniel A. Lazar}
	is a PhD student in the Department of Electrical and Computer Engineering at the University of California, Santa Barbara. In 2014, he received a B.Sc. in Electrical Engineering from Washington University in St. Louis, after which he spent two years working in wireless communications. His research focuses on control of transportation networks. 
\end{IEEEbiographynophoto}

\begin{IEEEbiographynophoto}{Ramtin Pedarsani}
	is an Assistant Professor in ECE Department at the
	University of California, Santa Barbara. He received the B.Sc. degree in
	electrical engineering from the University of Tehran, Tehran, Iran, in 2009,
	the M.Sc. degree in communication systems from the Swiss Federal Institute
	of Technology (EPFL), Lausanne, Switzerland, in 2011, and his Ph.D. from
	the University of California, Berkeley, in 2015. His research interests include
	networks, game theory, machine learning, and transportation
	systems. Ramtin is a recipient of the IEEE international conference on
	communications (ICC) best paper award in 2014. 
\end{IEEEbiographynophoto} 

\begin{IEEEbiographynophoto}{Sam Coogan}
	is an assistant professor at Georgia Tech with a joint appointment in the School of Electrical and Computer Engineering and the School of Civil and Environmental Engineering. Prior to joining Georgia Tech in 2017, he was an assistant professor in the Electrical Engineering Department at UCLA from 2015-2017. He received the B.S. degree in Electrical Engineering from Georgia Tech and the M.S. and Ph.D. degrees in Electrical Engineering from the University of California, Berkeley. His research is in the area of dynamical systems and control and focuses on developing scalable tools for verification and control of networked, cyber-physical systems with an emphasis on transportation systems. He received a CAREER award from the National Science Foundation in 2018, the Outstanding Paper Award for the IEEE Transactions on Control of Network Systems in 2017, and the Best Student Paper Award for the Hybrid Systems: Computation and Control conference in 2015.
\end{IEEEbiographynophoto}